\input vanilla.sty \scaletype{\magstep1}
\scalelinespacing{\magstep1} \def\bull{\vrule height .9ex
width .8ex depth -.1ex}


\title Calder\'on couples of re-arrangement invariant spaces
\endtitle

\author N.J.  Kalton\footnote{Research supported by
NSF-grants DMS-8901636 and DMS-9201357; the author also
acknowledges support from US-Israel BSF-grant 87-00244}\\
Department of Mathematics\\ University of
Missouri-Columbia\\ Columbia, Missouri 65211.  \endauthor

\subheading{Abstract}We examine conditions under which a
pair of re-arrangement invariant function spaces on $[0,1]$
or $[0,\infty)$ form a Calder\'on couple.  A very general
criterion is developed to determine whether such a pair is a
Calder\'on couple, with numerous applications.  We give, for
example, a complete classification of those spaces $X$ which
form a Calder\'on couple with $L_{\infty}.$ We specialize
our results to Orlicz spaces and are able to give necessary
and sufficient conditions on an Orlicz function $F$ so that
the pair $(L_F,L_{\infty})$ forms a Calder\'on pair.

\subheading{1.  Introduction}

Suppose $(X,Y)$ is a compatible pair of Banach spaces (see
[4] or [5]).  We denote by $K(t,f)=K(t,f;X,Y)$ the Peetre
K-functional on $X+Y$ i.e.  $$
K(t,f)=\inf\{\|x\|_X+t\|y\|_Y:\ x+y=f\}.$$ Then $(X,Y)$ is
called a {\it Calder\'on couple} (or a {\it
Calder\'on-Mityagin couple}) if whenever $f,g$ satisfy
$$K(t,f)\le K(t,g)$$ for all $t$ then there is a bounded
operator $T:X+Y\to X+Y$ such that $\|T\|_X,\|T\|_Y<\infty$
and $Tg=f.$ We will say that $(X,Y)$ is a {\it uniform}
Calder\'on couple (with constant $C$) if we can further
insist that $\max(\|T\|_X,\|T\|_Y)\le C.$ Calder\'on couples
are particularly important in interpolation theory because
it is possible to give a complete description of all
interpolation spaces for such a couple.  Indeed, for such a
couple, it is easy to show that a space $Z$ is an
interpolation space if and only if it is K-monotone, i.e. if
$f\in Z$ and $g\in X+Y$ with $K(t,g)\le K(t,f)$ then $g\in
Z.$ It follows from the K-divisibility theorem of Brudnyi
and Krugljak [7] that if $Z$ is a normed K-monotone space
then allows $\|f\|_Z$ on $Z$ is equivalent to a norm
$\|K(t,f)\|_{\Phi}$ where $\Phi$ is an appropriate lattice
norm on functions on $(0,\infty).$ Thus, for Calder\'on
couples, one has a complete description of all interpolation
spaces.  We remark also at this point that there are
apparently no known examples of Calder\'on couples which are
not uniform.

There has been a considerable amount of subsequent effort
devoted to classifying Calder\'on couples of
rearrangement-invariant spaces on $[0,1]$ or $[0,\infty).$
It is a classical result of Calder\'on and Mityagin
([9],[32]) that the pair $(L_1,L_{\infty})$ is a uniform
Calder\'on couple with constant 1. It is now known that any
pair $(L_p,L_q)$ is a Calder\'on couple (and indeed weighted
versions of these theorems are valid); we refer the reader
to Lorentz and Shimogaki [27], Sparr [36], Arazy and Cwikel
[1], Sedaev and Semenov [35] and Cwikel [13], [15].
Subsequent work has shown that under certain hypotheses
pairs of Lorentz spaces or Marcinkiewicz spaces are
Calder\'on couples; see Cwikel [14], Merucci [30], [31] and
Cwikel-Nilsson [16], [17].  For further positive results on
Calder\'on couples see [18] (for weighted Banach lattices),
and [21] and [38] (for Hardy spaces).

On the negative side, Ovchinnikov [34] showed that on
$[0,\infty)$ the pair $(L_1+L_{\infty},L_1\cap L_{\infty})$
is not a Calder\'on couple; indeed Maligranda and
Ovchinnikov show that if $p\neq 2$ then $L_p\cap L_q$ and
$L_p+L_q$ ($\frac1p+\frac1q=1$) are interpolation spaces not
obtainable by the K-method [29].

The general problem we consider in this article is that of
providing necessary and sufficient conditions on a pair
$(X,Y)$ of r.i. spaces (always assumed to have the so-called
Fatou property) on either $[0,1]$ or $[0,\infty)$ so that
$(X,Y)$ is a Calder\'on couple.  Although we cannot provide
a complete answer to this problem, we can resolve it in
certain cases and this enables us to settle some open
questions in the area (see, e.g.  Maligranda [28], Problems
1-3 or Brudnyi-Krugljak [8] p.685 [g], [i]).  For example,
we give a complete classification of all r.i. spaces $X$ so
that $(X,L_{\infty})$ is a Calder\'on couple and hence give
examples of r.i. spaces (even Orlicz spaces) $X$ so that
$(X,L_{\infty})$ is not a Calder\'on couple.  Our methods
give fairly precise information in the problem of
classifying pairs of Orlicz spaces which form Calder\'on
couples.  It should also be mentioned that our results apply
equally to symmetric sequence spaces.

We now describe our results in more detail.  Let $X$ be an
r.i. space on $[0,1]$ or $[0,\infty)$ or a symmetric
sequence space.  Let $e_n=\chi_{[2^{n},2^{n+1})}$ for
$n\in\bold J$ where $\bold J=\bold Z_-=-\bold N$ or $\bold
J=\bold Z$ or $\bold J=\bold N\cup \{0\}.$ We associate to
$X$ a K\"othe sequence space $E_X$ on $\bold J$ by defining
$$ \|\xi\|_{E_X} =\|\sum_{n\in\bold J}\xi(n)e_n\|_X.$$ We
then say that $X$ is {\it stretchable} if the sequence space
$E_X$ has the right-shift property (RSP) i.e. there is a
constant $C$ so that if $(x_n,y_n)_{n=1}^N$ is any pair of
finite normalized sequences in $E_X$ so that supp $x_1<$
supp $y_1<$ supp $x_2<\cdots<$ supp $y_n$ then for any
$\alpha_1,\ldots,\alpha_N$ we have $$
\|\sum_{n=1}^N\alpha_ny_n\|_{E_X} \le
C\|\sum_{n=1}^N\alpha_nx_n\|_{E_X}.$$ Thus $E_X$ has (RSP)
if the right-shift operator is uniformly bounded on the
closed linear span of every block basic sequence with
respect to the canonical basis.  We similarly say that $X$
is {\it compressible} if $E_X$ has the corresponding {\it
left-shift property} (LSP).  Finally we say that $X$ is {\it
elastic} if it is both stretchable and compressible.  It is
easy to see that $L_p-$spaces and more generally Lorentz
spaces with finite Boyd indices are elastic because $E_X$ in
this case is a weighted $\ell_p-$space (in fact this
property characterizes Lorentz spaces when the Boyd indices
are finite).  On the other hand it is not difficult to give
examples of r.i. spaces which are neither compressible nor
stretchable.  Curiously, however, we have no example of a
space which is either stretchable or compressible and not
elastic.

The significance of these ideas is illustrated by Theorem
5.4.  The pair $(X,L_{\infty})$ is a Calder\'on couple if
and only if $X$ is stretchable.  Dually if we assume that
$X$ has nontrivial concavity then $(X,L_1)$ is a Calder\'on
couple if and only if $X$ is compressible (Theorem 5.5).
More generally if $(X,Y)$ is any pair of r.i. spaces such
that either the Boyd indices satisfy $p_Y>q_X$ or there
exists $p$ so that $X$ is $p$-concave and $Y$ is $p$-convex
and has nontrivial concavity then $(X,Y)$ is a Calder\'on
couple if and only if $X$ is stretchable and $Y$ is
compressible.

In Section 6 we study these concepts for Orlicz spaces.  We
show that for an Orlicz space to be compressible it is
necessary and sufficent that it is stretchable; thus we need
only consider elastic Orlicz spaces.  We show for example
that $L_F[0,1]$ (where $F$ satisfies the
$\Delta_2-$condition) is elastic if and only if there is a
constant $C$ and a bounded monotone increasing function
$w(t)$ so that for any $0<x\le 1$ and any $1\le s\le
t<\infty$ we have $$ \frac{F(tx)}{F(t)} \le
C\frac{F(sx)}{F(s)} + w(t)-w(s).$$ This condition implies
that the Boyd indices (or Orlicz-Matuszewska indices) $p_F$
and $q_F$ of $L_F$ coincide.  In fact it implies the
stronger condition that $F$ must be equivalent to a function
which is regularly varying in the sense of Karamata (see
[6]).  We give examples to show that $F$ can be regularly
varying with $L_F$ inelastic and that $L_F$ can be elastic
without coinciding with a Lorentz space (cf.  [26], [33]).

Brudnyi (cf.  [8]) has conjectured that if a pair of
(distinct) Orlicz spaces $(L_F[0,1],\allowmathbreak
L_G[0,1])$ is a Calder\'on couple then $p_F=q_F$ and
$p_G=q_G.$ We show by example that this is false.  However
we also show that either $p_F=p_G$ and $q_F=q_G$ or both
$L_F$ and $L_G$ are elastic and hence $p_F=q_F$ and
$p_G=q_G.$

Let us now introduce some notation and conventions.  Let
$\Omega$ be a Polish space and let $\mu$ be a
$\sigma-$finite Borel measure on $\Omega.$ Let $L_0(\mu)$
denote the space of all real-valued Borel functions on
$\Omega$ (where functions differing on a set of measure zero
are identified), equipped with the topology of convergence
in $\mu-$measure on sets of finite measure.  By a K\"othe
function space on $\Omega$ we shall mean a Banach space $X$
which is a subspace of $L_0$ containing the characteristic
function $\chi_B$ whenever $\mu(B)<\infty$ and such that the
norm $\|\,\|_X$ satisfies the conditions:\newline (a)
$\|f\|_X \le \|g\|_X$ whenever $|f|\le |g|$ a.e.\newline (b)
$B_X=\{f:\|f\|_X\le 1\}$ is closed in $L_0.$\newline
Condition (b) is usually called the Fatou property; note
here that we include the Fatou property in our definition
and so it is an implicit assumption throughout the paper.
It is sometimes convenient to extend the definition of
$\|f\|_X$ by setting $\|f\|_X=\infty$ if $f\notin X.$ We
will also write $P_Bf=f_B=f\chi_B$ when $B$ is a Borel
subset of $\Omega.$ We let supp $f=\{\omega:f(\omega)\neq
0\}.$

If $X$ is a K\"othe function space then we say that $X$ is
$p$-convex ($1\le p\le \infty)$ if there is a constant $M$
so that for any $f_1,\ldots f_n\in X$ we have $$
\|(\sum_{k=1}^n |f_k|^p)^{1/p}\|_X \le M(\sum_{k=1}^n
\|f_k\|_X^p)^{1/p}$$ and $p$-concave if there exists $M$ so
that $$ (\sum_{k=1}^n\|f_k\|_X^p)^{1/p} \le M
\|(\sum_{k=1}^n|f_k|^p)^{1/p}\|_X.$$ Similarly we say that
$X$ has an upper $p$-estimate if there is a constant $M$ so
that if $f_1,\ldots,f_n$ are disjoint in $X$ then $$
\|\sum_{k=1}^nf_k\|_X \le M(\sum_{k=1}^n\|f_k\|_X^p)^{1/p}$$
and $X$ has a lower $p$-estimate if there exists $M$ so that
if $f_1,\ldots,f_n$ are disjoint then $$
(\sum_{k=1}^n\|f_k\|_X^p)^{1/p}\le M\|\sum_{k=1}^nf_k\|_X.$$
See Lindenstrauss-Tzafriri [25] for a fuller discussion.

We will sometimes use $\langle f,g\rangle$ for
$\int_{\Omega}fg\,d\mu.$ With this notion of pairing we will
also use $X^*$ for the K\"othe dual of $X$ (which will
coincide with the full dual if $X$ is separable).

If $(X,Y)$ are two K\"othe function spaces on $(\Omega,\mu)$
then the pair $(X,Y)$ is necessarily Gagliardo complete (cf.
[4]).  We denote by $\Cal A(X,Y)$ the space of admissible
operators i.e. operators $T:X+Y\to X+Y$ such that
$\|T\|_X=\sup\{\|Tf\|_X:\|f\|_X\le 1\}<\infty$ and
$\|T\|_Y=\sup\{\|Tf\|_Y:\|f\|_Y\le 1\}<\infty.$ We norm
$\Cal A(X,Y)$ by $\|T\|_{(X,Y)}=\max(\|T\|_X,\|T\|_Y).$

In the special case when $\Omega=\bold J$ is a subset of
$\bold Z$ and $\mu$ is counting measure we write
$\omega(\bold J)=L_0(\mu)$ and a K\"othe function space $X$
is called a K\"othe sequence space modelled on $\bold J.$ An
operator $T$ on $X$ is then called a matrix if it takes the
form $$ Tx(n) =\sum_{k\in\bold J}a_{nk}x(k)$$ for some
$(a_{nk})_{n,k\in\bold J}.$ We remark here that the
assumption that $T$ is a matrix forces the existence of an
adjoint operator $T^*:X^*\to X^*$ even in the nonseparable
situation when $X^*$ is not the full dual of $X.$

If $\Omega=[0,1]$ or $[0,\infty)$ (with $\mu$ Lebesgue
measure) or if $\Omega=\bold N$ (with $\mu$ counting
measure) then for $f\in L_0(\mu)$ we define the decreasing
rearrangement $f^*$ of $f$ by $$f^*(t) = \sup_{B:\mu(B)\le
t}\inf_{s\in B}|f(s)|,$$ for $0< t<\infty.$ We say that $X$
is a rearrangement-invariant space (or a symmetric sequence
space if $\Omega=\bold N)$ if $\|f\|_X=\|f^*\|_X$ for all
$f\in L_0.$ If we define $$f^{**}(t)
=\frac1t\int_0^tf^*(s)ds$$ then it is well-known that if
$f,g\in L_0$ with $f^{**}\le g^{**}$ then $\|f\|_X\le
\|g\|_X.$

If $X$ is an r.i. space on $[0,1]$ or $[0,\infty)$ then the
dilation operators $D_a$ on $X$ are then defined by $
D_af(t) = f(t/a)$ (where we regard $f$ as vanishing outside
$[0,1]$ in the former case).  We can then define the Boyd
indices $p_X$ and $q_X$ of $X$ by $$ \align p_X
&=\lim_{a\to\infty} \frac{\log a}{\log \|D_a\|_X} \\ q_X
&=\lim_{a\to 0}\frac{\log a}{\log \|D_a\|_X}.  \endalign $$
In the case when $X$ is a symmetric sequence space we define
$p_X$ and $q_X$ in the same way but we define $D_a$ by the
nonlinear formula $$ D_af(n) = f^*(n/a)$$ where $f^*$ is
well-defined on $[0,\infty).$

Finally let us mention two special classes of r.i. spaces.
If $1\le p< \infty$ we will say that an r.i. space $X$ on
$\Omega=[0,1]$ or $[0,\infty)$ is a {\it Lorentz space of
order $p$} if there is a positive monotone increasing weight
function $w:\Omega\to (0,\infty)$ such that
$\sup_{t,2t\in\Omega}w(2t)/w(t)<\infty$ and $\|f\|_X$ is
equivalent to the quasinorm $$ \|f\|_{w,p} =
(\int_{\Omega}f^*(t)^pw(t)^p\frac{dt}t)^{1/p}.$$ We can then
write $X=L_{w,p}.$ If we take $w(t)=t^{1/q}$ we obtain the
standard Lorentz spaces $L(q,p).$ It is easy to compute that
the Lorentz space $X=L_{w,p}$ has Boyd indices $p_X,q_X$
where $$ \align \frac1{p_X} &=
\lim_{a\to\infty}\sup_{t,at\in\Omega} \frac{\log w(at)-\log
w(t)}{\log a}\\ \frac1{q_X} &= \lim_{a\to\infty}
\inf_{t,at\in\Omega}\frac{\log w(at)-\log w(t)}{\log a}.
\endalign $$ If we impose the additional restriction that
$q_X<\infty$ then it can easily be seen that we may suppose
that $w$ satisfies $\inf_{t,2t\in\Omega}w(2t)/w(t)>1.$

We will also be interested in Orlicz function spaces and
sequence spaces.  By an Orlicz function we shall mean a
continuous strictly increasing convex function
$F:[0,\infty)\to[0,\infty)$ such that $F(0)=0.$ $F$ is said
to satisfy the $\Delta_2-$condition if there is a constant
$\Delta$ such that $F(2x)\le \Delta F(x)$ for all $x\ge 0.$

The Orlicz function space $L_F(\Omega,\mu)$ is defined by $$
\|f\|_{L_F}
=\inf\{\alpha>0:\int_{\Omega}F(\alpha^{-1}f(t))dt\le 1 \}$$
so that $L_F=\{f:\|f\|_{L_F}<\infty\}.$

In this case the Boyd indices $p_F=p_{L_F}$ and
$q_F=q_{L_F}$ are closely related to the Orlicz-Matuszewska
indices of $F$ (see Lindenstrauss-Tzafriri [25] p. 139).
More precisely let $\alpha^{\infty}(F)$ (resp.
$\alpha^0(F)$) be the supremum of all $p$ so that for some
$C$ we have $F(st)\le Cs^pF(t)$ for all $0\le s\le 1$ and
all $t\ge 1$ (resp.  $t\le 1$).  Similarly let
$\beta^{\infty}(F)$ (resp.  $\beta^{0}(F)$) be the infimum
of all $q$ so that for some $C$ we have $s^pF(t) \le CF(st)$
for all $0\le s\le 1$ and all $t\ge 1$ (resp.  $t\le 1$).
Then if $\Omega=[0,1]$ we have $p_F=\alpha^{\infty}(F)$ and
$q_F=\beta^{\infty}(F).$ If $\Omega=[0,\infty)$ then
$p_F=\min(\alpha^{\infty}(F),\alpha^0(F))$ and
$q_F=\max(\beta^{\infty}(F),\beta^0(F)).$ If we assume the
$\Delta_2-$condition (and we always will) then $q_F<\infty.$

\vskip2truecm

\subheading{2.  The shift properties}

Let $\bold J$ be one of the three sets $\bold Z$, $\bold
Z_+=\{n\in\bold Z:  n\ge 0\}$ or $\bold Z_{-}=\bold
Z\setminus \bold Z_+.$ Let $\omega(\bold J)$ denote the
space of all sequences modelled on $\bold J.$ If
$x=\{x(k)\}_{k\in\bold J}$ is a sequence (modelled on $\bold
J$) we write supp $x=\{k:x(k)\neq 0\}.$ If $A,B$ are subsets
of $\bold J$ we write $A<B$ if $a<b$ for every $a\in A,b\in
B.$ If $I$ is any interval of $\bold Z$ and $(x_n,y_n)_{n\in
I}$ is a pair of sequences in $\omega(\bold J)$ we say
$(x_n,y_n)$ is {\it interlaced} if each $x_n,y_n$ has finite
support and supp $x_n <$ supp $y_n\ (n\in I)$ and supp
$y_n<$ supp $x_{n+1}$ whenever $n,n+1\in I.$

Let $E$ be a K\"othe sequence space modelled on $\bold J.$
We will say that $E$ has the {\it right-shift property
(RSP)} if there is a constant $C$ such that whenever
$(x_n,y_n)_{n\in I}$ is an interlaced pair with $\|y_n\|_E
\le \|x_n\|_E=1\ (n\in I)$ then for every finitely nonzero
sequence of scalars $(\alpha_n)_{n\in I}$ we have $$
\|\sum_{n\in I} \alpha_ny_n\|_E \le C\|\sum_{n\in
I}\alpha_nx_n\|_E.$$

Conversely we will say that $E$ has the {\it left-shift
property (LSP)} if there is a constant $C'$ so that for
every interlaced pair $(x_n,y_n)_{n\in I}$ with
$\|x_n\|_E\le \|y_n\|_E=1,$ and every finitely nonzero
$(\alpha_n)_{n\in I}$ we have $$ \|\sum_{n\in I
}\alpha_nx_n\|_E \le C'\|\sum_{n\in I}\alpha_ny_n\|_E.$$

\proclaim{Proposition 2.1}$E$ has (LSP) if and only if $E^*$
has (RSP).\endproclaim

\demo{Proof}We will only prove one direction.  Let us assume
$E^*$ has (RSP) with constant $C$.  Let $(x_n,y_n)_{n\in
I},$ be an interlaced pair with $\|y_n\|_E\le \|x_n\|_E=1.$
We may assume each $x_n,y_n$ is positive (i.e.
$x_n(k),y_n(k)\ge 0$ for every $k$).  Suppose
$(\alpha_n)_{n\in I}$ is a finitely nonzero sequence of
nonnegative reals.  Let $f=\sum\alpha_ny_n$.  Then there
exists positive $g\in E^*$ with supp $g\subset$supp $f$ and
so that $\langle f,g\rangle=\|f\|_E$ while $\|g\|_{E^*}=1$.
We can write $$ g =\sum \beta_n v_n$$ where each $v_n$ is
positive, $\|v_n\|_{E^*}=1$ and supp $v_n\subset$ supp
$y_n.$

Next pick positive $u_n$ with supp $u_n\subset$ supp $x_n$,
$\langle x_n,u_n\rangle=1$ and $\|u_n\|_{E^*}=1.$ We
conclude from the fact that $E^*$ has (RSP) that $$
\|\sum_{n\in I}\beta_n u_n\|_{E^*}\le C.$$ Thus $$ \align
\|\sum_{n\in I}\alpha_ny_n\|_E &= \sum_{n\in I
}\alpha_n\beta_n\langle y_n,v_n\rangle\\ &\le \sum
\alpha_n\beta_n \\ &\le \langle
(\sum\alpha_nx_n),(\sum\beta_nu_n)\rangle\\ &\le
C\|\sum\alpha_nx_n\|_E.  \endalign $$ Thus the proposition
is proved.\bull\enddemo

\proclaim{Proposition 2.2}Suppose $E$ is a K\"othe sequence
space modelled on $\bold Z.$ Define $E_+=E(\bold Z_+)$ and
$E_-=E(\bold Z_-).$ Then $E$ has (RSP) (resp.  (LSP)) if and
only both $E_+$ and $E_-$ have (RSP) (resp.
(LSP)).\endproclaim

\demo{Proof}One direction is obvious.  For the other,
suppose both $E_+$ and $E_-$ have (RSP) with constant $C$,
say.  Suppose $(x_n,y_n)_{n\in I}$ is an interlaced pair of
sequences with $\|y_n\|_E\le \|x_n\|_E=1$ and that
$(\alpha_n)_{n\in I}$ is finitely nonzero.  Then there
exists $m\in I$ so that supp $(x_n+y_n)\subset \bold Z_-$
for $n<m$ and supp $(x_n+y_n)\subset \bold Z_+$ for $n>m.$
Now $$ \align \|\sum_{n\in I}\alpha_ny_n\|_E &\le |\alpha_m|
+\|\sum_{n<m}\alpha_ny_n\|_E + \|\sum_{n>m}\alpha_ny_n\|_E\\
&\le (2C+1) \|\sum_{n\in I}\alpha_nx_n\|_E.  \endalign $$
Thus $E$ has (RSP) with constant at most
$2C+1.\bull$\enddemo

To simplify our discussion we introduce the idea of an
order-reversal.  Let $E=E(\bold J)$ be a K\"othe sequence
space.  We let $\tilde{\bold J}=\{-(n+1):n\in \bold J\}$ and
if $x\in\omega(\bold J)$ we set $\tilde x(n)=x(-(n+1))$ for
$n\in\bold J.$ Let $\tilde E(\tilde{\bold J})$ be defined by
$\|x\|_{\tilde E}=\|\tilde x\|_E$; then $\tilde E$ is the
order-reversal of $E.$ Clearly (LSP) (resp.  (RSP)) for $E$
is equivalent to (RSP) (resp.  (LSP)) for $\tilde E.$

Next observe that if $(w_n)_{n\in \bold J}$ satisfy $w_n>0$
for all $n$ then the weighted sequence space $E(w)=\{x:xw\in
E\}$ normed by $\|x\|_{E(w)}=\|xw\|_E$ satisfies (LSP)
(resp.  (RSP)) if and only if $E$ satisfies (LSP) (resp.
(RSP)).

\proclaim{Proposition 2.3}Let $E=E(\bold J)$ be a symmetric
sequence space.  Suppose $E$ has either (LSP) or (RSP).
Then $E=\ell_p(\bold J)$ for some $1\le
p\le\infty.$\endproclaim

\demo{Proof}For convenience of notation we consider only the
case $\bold J=\bold Z_+$ and (LSP) and leave the reader to
make the minor adjustments necessary for the other cases.
Let $(u_n)_{n\in\bold N}$ be any normalized positive block
basic sequence in $E(\bold J)$.  Select $a_n\in$ supp $u_n.$
Then $(u_{2n},e_{a_{2n+1}})_{n \in \bold N}$ is an
interlaced pair.  Thus $$ \|\sum_{n\in\bold N}
\alpha_ne_n\|_E \le C\|\sum_{n\in\bold
N}\alpha_nu_{2n}\|_E.$$ Similarly
$(e_{a_{2n-1}},u_{2n})_{n\in\bold N}$ is an interlaced pair
and so $$ \|\sum_{n\in\bold N}\alpha_n u_{2n}\|_E \le
C\|\sum_{n\in\bold N}\alpha_ne_n\|_E.$$ Thus $(u_{2n})$ is
$C^2-$equivalent to $(e_n)$ and similarly so is
$(u_{2n-1})_{n\in\bold N}$.  It then follows that the basis
(or basic sequence) $(e_n)$ is perfectly homogeneous and by
a theorem of Zippin [39] (see Lindenstrauss-Tzafriri [24])
this implies that it is equivalent to the $\ell_p$-basis for
some $p$ or the $c_0$-basis; in the latter case we deduce
that $E=\ell_{\infty}(\bold J).$ The result then follows.
\bull\enddemo

\proclaim{Proposition 2.4}Let $E=E(\bold Z_+)$ be a K\"othe
sequence space with (LSP) or (RSP).  If $E$ contains a
symmetric basic sequence then there exists $1\le p\le\infty$
and an increasing sequence $(a_k)_{k\ge 0}$ with $a_0=0$ so
that $E=\ell_p(E[a_k,a_{k+1})).$ In particular, when $E$ is
separable, we have $p<\infty$ and any symmetric basic
sequence in $E$ is equivalent to the canonical
$\ell_p-$basis.\endproclaim

\demo{Remark}Of course there is a similar result if $\bold
J=\bold Z_-.$ However in the two-ended setting $\bold
J=\bold Z$ we recall that $E$ has (RSP) (resp.  (LSP)) if
and only if both $E(\bold Z_+)$ and $E(\bold Z_-)$ have
(RSP) (resp.  (LSP)).  In particular, $\ell_p(\bold
Z_-)\oplus\ell_r(\bold Z_+)$ has (LSP) and (RSP) even if
$r\neq p.$\enddemo

\demo{Proof}We can suppose that $(u_n)$ is a normalized
symmetric block basic sequence.  By an interlacing argument
as in Proposition 2.3 it will follow that a subsequence
$(e_{a_k})$ of the unit vectors is symmetric.  Since the
restriction of $E$ to this subsequence has (LSP) or (RSP) it
follows that it is equivalent to the $\ell_p$-basis for some
$1\le p<\infty$ or to the $c_0-$basis by Proposition 2.3.
For convenience we suppose the former case and fix $p.$ Let
$I_k= I_k=[a_k,a_{k+1})$.  Then, for suitable $C$, by an
interlacing argument any normalized sequence $(v_k)$
supported on $I_{2k}$ is $C-$equivalent to the $\ell_p-
$basis; similarly any normalized sequence supported on
$I_{2k+1}$ is $C-$equivalent to the $\ell_p-$basis and the
first part of the result follows.  For the last part, if $E$
is separable then obviously $p<\infty$ and a simple blocking
argument gives the result.  \bull\enddemo

\demo{Remark}It is possible that $E$ contains no symmetric
basic sequence.  Indeed, Tsirelson space $T$ [37] and its
convexifications provide examples of such spaces with (RSP)
and (LSP) (see [10] and [12]).  It is not difficult to use
Krivine's theorem [22] to show that if $E=E(\bold Z_+)$ has
(LSP) (or (RSP)) then there is a subsequence $(e_{a_n})$ of
the unit vector basis so that for some $C,p$ we have for all
$k$ and every $k$ vectors $x_1,x_2,\ldots,x_k$ with supp
$x_1<$ supp $x_2<\ldots<$ supp $x_k$ and supp
$(x_1+\cdots+x_k)\subset \{a_n\}_{n\ge 1}$ then $$
C^{-1}(\sum_{n=1}^k\|x_n\|_E^p)^{1/p} \le
\|\sum_{n=1}^kx_n\|_E\le
C(\sum_{n=1}^{\infty}\|x_n\|_E^p)^{1/p}$$ with appropriate
modifications when $p=\infty.$ Thus any space $E$ having
either (LSP) or (RSP) and no symmetric basic sequence has a
``Tsirelson-like'' subspace.

\demo{Problem}Does there exist a K\"othe sequence space with
(LSP) and not (RSP)?

Let us remark that this is probably non-trivial.  Indeed the
corresponding question for simple shifts has been considered
[3] and Bellenot has only recently given an example
[2].\enddemo

\proclaim{Lemma 2.5}Let $E$ be a K\"othe sequence space on
$\bold J$ with (RSP); then there is a constant $C$ so that
whenever $(x_n,y_n)_{n\in I}$ is an interlaced pair of
sequences with $\|y_n\|_E\le\|x_n\|_E=1$ and $(x_n^*)_{n\in
I}$ is a sequence in $E^*$ with supp $x_n^*\subset$ supp
$x_n$ and $x_n^*(x_n)=\|x_n\|_{E^*}=1$ then the operator $T$
defined by $Tx =\sum_{n\in I}\langle x,x_n^*\rangle y_n$ is
bounded on $E$ with $\|T\|_E\le C.$\endproclaim

\demo{Proof}For any $x\in E$ with finite support, $Tx$ has
finite support and we can define $g_n^*\in E^*$ and a
finitely nonzero sequence $(\alpha_n)_{n\in I}$ so that
$\|g_n^*\|=1 (n\in I)$, supp $g_n^*\subset $ supp $y_n$,
$\|\sum \alpha_ng_n^*\|_{E^*}=1$ and $$ \langle
Tx,\sum_{n\in I}\alpha_ng_n^*\rangle =\|Tx\|_E.$$ Thus $$
\align \|Tx\|_E &= \sum_{n\in I}\alpha_nx_n^*(x)\\ &=
\langle x, \sum_{n\in I}\alpha_n x_n^*\rangle \\ &\le
\|x\|_E \|\sum_{n\in I}\alpha_nx_n^*\|_{E^*}\\ &\le C\|x\|_E
\|\sum_{n\in I}\alpha_ng_n^*\|_{E^*}\\ &\le C\|x\|_E
\endalign $$ where $C$ is the (LSP) constant of $E^*$ (which
actually is the (RSP) constant of $E$ by Proposition 2.1 and
its proof).  The result follows.\bull\enddemo

\proclaim{Lemma 2.6}Under the hypotheses of Lemma 2.5, there
is a constant $C_1$ so that $(J_n)_{n\in I}$ is a sequence
of intervals in $\bold J$ with $J_n<J_{n+1}$ whenever
$n,n+1\in I$, $(x_n)_{n\in I},(y_n)_{n\in I}$ are two
normalized sequences with supp $x_n,$ supp $y_n\subset J_n$
and $(x_n^*)$ is any sequence with supp $x_n^*\subset$ supp
$x_n$ and $x_n^*(x_n)=1=\|x_n^*\|_{E^*}$ then the operator
$$ Tx =\sum_{n\in I}\langle x,x_n^*\rangle y_{n+1}$$ (where
$y_{n+1}=0$ if $n+1\notin I$) is bounded on $E$ with
$\|T\|_E\le C_1.$\endproclaim

\demo{Proof}The sequence pairs
$(x_{2n},y_{2n+1})_{2n,2n+1\in I}$ and
$(x_{2n-1},y_{2n})_{2n-1,2n\in I}$ are interlaced and the
lemma follows from 2.5 with $C_1=2C$ by simply
adding.\bull\enddemo

\demo{Remark}If $E$ is separable and has both (LSP) and
(RSP) then Lemma 2.6 quickly shows that every normalized
block basic sequence in $E$ spans a complemented subspace;
this property is, of course, enjoyed by Tsirelson space [12]
(see also Casazza-Lin [11] for an earlier similar example).
If this property holds for a symmetric sequence space then
it is isomorphic to $\ell_p$ for some $1\le p<\infty$ (see
Lindenstrauss-Tzafriri [23]).  \enddemo

\vskip2truecm

\subheading{3.  The shift properties for pairs}

We next consider a pair of K\"othe sequences spaces $(E,F)$
modelled on $\bold J.$ We will say that $(E,F)$ has the
(RSP) if there is a constant $C$ so that whenever
$\{x_n,y_n\}_{n\in I}$ is an interlaced pair with $\|y_n\|_E
\le \|x_n\|_E=1$ and $x_n,y_n\ge 0$ then there is a positive
matrix $T$ with $\|T\|_{(E,F)}\le C$ and $Tx_n=y_n.$ We will
say that the pair $(E,F)$ has (LSP) if $(\tilde F,\tilde E)$
has (RSP).  If $(E,F)$ has both (LSP) and (RSP) then we say
that it has the {\it shift property} (SP).

We first note that if $(E,F)$ has (RSP) then $E$ has (RSP).
Conversely it follows from Lemma 2.5 that if $E$ has (RSP)
then $(E,E)$ has (RSP).

In this section, we show that, under certain hypotheses, one
can deduce (RSP) for the couple $(E,F)$ from the property
(RSP) for E alone.  We will need some definitions.  We
define the shift operators $\tau_n$ for $n\in\bold Z$ on
$\omega(\bold J)$ by $\tau_n(x)(k)=x(k-n),$ where we
interpret $x(j)=0$ when $j\notin\bold J.$ We define
$\kappa_+(E)=\lim_{n\to\infty}\|\tau_n\|_E^{1/n}$ (which can
be $\infty$ in the case when $\tau_1$ is unbounded on $E$)
and $\kappa_-(E)=\lim_{n\to\infty}\|\tau_{-n}\|_E^{1/n}.$ We
will also let $\rho(n)=\rho(n;E,F)=\|e_n\|_E/\|e_n\|_F.$ We
will say that $(E,F)$ is {\it exponentially separated} if
there exists $\beta>0$ and $C_0$ so that if $m,m+n\in\bold
J$ then $\rho(m+n)\ge C_0^{-1}2^{n\beta }\rho(m).$

\proclaim{Lemma 3.1}If $\kappa_-(E)\kappa_+(F)<1$ then
$(E,F)$ is exponentially separated.\endproclaim

\demo{Proof}Here we have $\rho(m+n)/\rho(m)\ge
(\|\tau_{-n}\|_{E }\|\tau_n\|_{F})^{-1}.$ The hypothesis
then implies that for some $\beta>0$ we have
$\|\tau_{-n}\|_{E}\|\tau_n\|_{F}\le C2^{-n\beta}$ for some
$C>0.$ The result then follows.\bull\enddemo

\proclaim{Lemma 3.2}Suppose $(E,F)$ is exponentially
separated.  Then there is a constant $C_1$ so that if supp
$x\subset [a,b]$ then $$ C_1^{-1}\rho(a)\|x\|_F \le \|x\|_E
\le C_1\rho(b)\|x\|_F.$$ \endproclaim

\demo{Proof}Suppose $x=\sum_{k=a}^bx(k)e_k.$ Then $$ \align
\|x\|_E &\le \sum |x(k)|\|e_k\|_E \\ &\le \sum
|x(k)|\rho(k)\|e_k\|_F\\ &\le C_0\sum
\rho(b)2^{-\beta(b-k)}|x(k)|\|e_k\|_F\\ &\le
C_0\rho(b)(\sum_{k=0}^{\infty}2^{-\beta k})\max_j
|x(j)|\|e_j\|_F\\ &\le C\rho(b)\|x\|_F \endalign $$ for a
suitable constant $C.$ The other inequality is
similar.\bull\enddemo

\proclaim{Lemma 3.3}Let $E,F$ be a pair of K\"othe sequence
spaces satisfying $\kappa_-(E)\kappa_+(F)<1.$ Suppose $E$
has (RSP).  Then $(E,F)$ has the (RSP).  Similarly if $F$
has (LSP) then $(E,F)$ has (LSP).\endproclaim

\demo{Proof}We first note that it is only necessary to prove
the first statement since $(\tilde F,\tilde E)$ will satisfy
the same hypothesis $\kappa_-(\tilde F)\kappa_+(\tilde E)<1$
and $\tilde F$ will have (RSP) if $F$ has (LSP).

We refer back to Lemma 2.5; it is clear that there exists
$C_0$ so that if $\{x_n,y_n\}_{n\in I}$ is a positive
interlaced pair with $\|y_n\|_E\le \|x_n\|_E=1$ then if we
pick $x_n^*\ge 0$ with supp $x_n^*\subset$ supp $x_n$ and
$\langle x_n,x_n^*\rangle=\|x_n^*\|_{E^*}=1$ for $n\in I$
then $\|T\|_E\le C_0$ where $$ Tx =\sum_{n\in A}\langle
x,x_n^*\rangle y_{n}.$$ Obviously $T$ is a positive matrix.
We now compute $\|T\|_F.$ Suppose $k\in \text{supp }y_n$
where $n\in I.$ Then, since supp $x_n^*\subset (-\infty,k)$
and $y_n(k)e_k\le y_n,$ $$ |Tx(k)| = |x_n^*(x)|y_{n}(k)\le
\|x_{(-\infty,k)}\|_E\|e_k\|_E^{-1}.$$ Now $$ \align
\|x_{(-\infty,k)}\|_E &\le \sum_{j<k}|x(j)|\|e_j\|_E \\ &\le
\|e_k\|_E\sum_{j<k}|x(j)|\|\tau_{j-k}\|_E.  \endalign $$ We
have $$ |Tx(k)| \le
\sum_{j<k}|\tau_{k-j}x(k)|\|\tau_{j-k}\|_E$$ and hence,
since $Tx(k)$ vanishes for all coordinates not of this form,
$$ |Tx| \le \sum_{j=1}^{\infty} \|\tau_{-j}\|_E|\tau^jx|.$$
The hypothesis $\kappa_-(E)\kappa_+(F)<1$ implies that there
exists $M$ and $0<\delta<1$ so that
$\|\tau_{-j}\|_E\|\tau_{j}\|_F \le M\delta^j$ for $j>0.$
Hence:  $$ \|Tx\|_F \le
\sum_{j=1}^{\infty}M\delta^{j}\|x\|_F$$ so that $\|T\|_F\le
C_1$ for some constant $C_1$ depending only on $E,F.$
\bull\enddemo

Although Lemma 3.3 is enough for most of our purposes, there
are some possible modifications.  First we give a simple
argument in the case $F=\ell_{\infty},$ which will be useful
later.

\proclaim{Lemma 3.4}Suppose $E$ is a K\"othe sequence space
with (RSP) and that $F=\ell_{\infty}(\bold J)$.  Suppose
$(E,F)$ is exponentially separated.  Then $(E,F)$ has (RSP).
\endproclaim

\demo{Proof}We may assume that for some $C_0,\beta>0$ we
have $\|e_m\|_E\le C_02^{-\beta n}\|e_{m+n}\|_E.$ In this
case we proceed as in Lemma 3.3 but note that $$ \|Tx\|_F
=\sup_{k\in J}|Tx(k)|.$$ If $k\in \text{supp }y_n,$ $$
\align |Tx(k)| &=|\langle x,x_n^*\rangle|y_{n}(k)\\ &\le
\|x_{(-\infty,k)}\|_E\|e_k\|_E^{-1} \\ &\le
\sum_{j<k}|x_j|\|e_j\|_E\|e_k\|_E^{-1}\\ &\le C_0\sum_{j<k}
|x_j|2^{-\beta(k-j)}\\ &\le C_1\|x\|_F \endalign $$ for a
suitable $C_1.$ \bull\enddemo

Another version of Lemma 3.3, which actually generalizes
Lemma 3.4, is given by:

\proclaim{Lemma 3.5}Suppose $(E,F)$ is exponentially
separated, $E$ has (RSP) and that either (a) there exists
$1\le p\le\infty$ so that $E$ has a lower $p$-estimate and
$F$ has an upper $p$-estimate or (b) $E$ is $r$-concave for
some $r<\infty$ and there exists $1<q< \infty$ so that $E$
has an upper $q$-estimate and $F$ has a lower $q$-estimate.
Then $(E,F)$ has (RSP).  \endproclaim

\demo{Proof} For (a) we note that the case $p=\infty$ is
essentially covered in Lemma 3.4.  Suppose $p<\infty.$ By
Lemma 3.2 there is a constant $C_1$ so that if supp $x<$
supp $y$ then $\|x\|_E\|y\|_F \le C_1\|x\|_F\|y\|_E.$ There
is also a constant $C_2$ so that if $u_1,\ldots,u_n$ are
disjoint vectors in $E$ or $F$, $$ \align (\sum_{j=1}^n
\|u_j\|_E^p)^{1/p} &\le C_2\|\sum_{j=1}^nu_j\|_E \\ \|
\sum_{j=1}^nu_j\|_F &\le C_2(\sum_{j=1}^n\|u_j\|_F^p)^{1/p}.
\endalign $$ We suppose $(x_n,y_n)_{n\in I}$ is a positive
interlaced pair with $\|y_n\|_E \le \|x_n\|_E=1.$ Define $T$
as in Lemma 3.3, and set $J_k=$ supp $x_k^*.$ Now if $x\in
F,$ $$ \align \|Tx\|_F &= \|\sum_{k\in
A}x_k^*(x)y_{k+1}\|_F\\ &\le C_2 (\sum_{k\in
A}\|x_{J_k}\|_E^p\|y_{k+1}\|_F^p)^{1/p}\\ &\le C_2C_1
(\sum_{k\in A} \|x_{J_k}\|_F^p)^{1/p}\\ &\le C_2^2C_1
\|x\|_F \endalign $$ and (a) follows.

We turn to the proof of (b).  Let $1<p<\infty$ be conjugate
to $q.$ Then $\tilde E^*$ has (RSP) by Proposition 2.1.
Also $(\tilde E^*,\tilde F^*)$ is exponentially separated,
and $\tilde E^*$ has a lower $p$-estimate while $\tilde F^*$
has an upper $p$-estimate; thus by (a) the couple $(\tilde
E^*,\tilde F^*)$ has (RSP).  It follows that $(F^*,E^*)$ has
(LSP).  We further can assume, by renorming, that $E$ has an
upper $p$-estimate with constant $1$ and an $r$-concavity
constant $1$ (apply Lindenstrauss-Tzafriri [25] p.88 Lemma
1.f.11 to $E^*$.)  Let $C_1$ be the associated (LSP)
constant for this couple.  We first prove a claim:

\proclaim{Claim} There exist constants $C_2$ and $\delta<1$
depending only on $(E,F)$ with the following property.
Suppose $\{x_n,y_n\}_{n\in I}$ is a positive interlaced pair
of sequences with $\|x_n\|_E=\|y_n\|_E=1$.  Then there is a
subset $D$ of $\bold J$, and a positive matrix operator $S$
with $\|S\|_{(E,F)}\le C_2$ so that $Sx_n=P_Dy_{n}$ and
$\|y_{n}-P_Dy_{n}\|_E\le \delta,$ whenever $n\in I.$
\endproclaim

Choose $x_k^*\ge 0$ with supp $x_k^*\subset$ supp $x_k$ and
$\|x_k^*\|_E=\|x_k\|_E=\langle x_k,x_k^*\rangle =1.$
Similarly choose $y_k^*\ge 0$ with supp $y_k^*\subset$ supp
$y_k$ and $\|y_k^*\|_E=\|y_k\|_E=\langle y_k,y_k^*\rangle.$
We begin by using the (LSP) property of $(F^*,E^*)$ to
produce a positive matrix $V$ on $E$ with $\|V\|_{(E,F)}\le
C_1$ and $V^*y_{k}^*=x_k^*$ whenever $k\in I.$ Thus $\langle
Vx_k,y_{k}^*\rangle=1.$

Fix $\tau>0$ small enough so that
$\frac12\tau-\frac1pC_1^p\tau^p=\gamma>0.$ Let $D_k$ be the
set of $j\in \text{supp }y_k$ so that $2Vx_k(j)\ge
y_{k}(j).$ Let $D=\cup_{k\in I} D_k$.  Clearly there is a
positive matrix $S$ with $\|S\|_{(E,F)}\le 2C_1=C_2$ and
$Sx_k=P_Dy_{k}.$ Now observe that $\langle
Vx_k-P_DVx_k,y_{k}^*\rangle\le \frac12$ so that $$ \langle
\tau P_DVx_k + y_{k}-P_Dy_{k},y_{k}^*\rangle\ge
1+\frac12\tau-\|P_Dy_{k}\|_E.$$ Thus $$
1+\frac12\tau-\|P_Dy_{k}\|_E \le (C_1^p\tau^p + 1)^{1/p}\le
1+\frac1p C_1^p\tau^p.$$ Upon reorganization this yields:
$$ \|P_Dy_{k}\|_E \ge \frac12\tau -\frac1p C_1^p\tau^p
=\gamma.$$ This in turn implies $$ \|y_{k}-P_Dy_{k}\|_E \le
(1-\gamma^r)^{1/r}=\delta<1.$$ This establishes the claim.

To complete the proof from the claim is quite easy by an
inductive argument.  We may clearly construct a disjoint
sequence of subsets $(D_n)_{n\ge 1}$ of $\bold J$ and a
sequence of positive matrix operators $(S_n)_{n\ge 1}$ with
$\|S_n\|_{(E,F)}\le 2C_1\delta^{n-1}$ and so that
$S_nx_k=P_{D_{n}}y_{k}$ and
$\|y_{k}-\sum_{j=1}^nP_{D_j}y_{k}\|_E \le \delta^n.$ The
operator $T=\sum_{n=1}^{\infty}S_n$ is a positive matrix and
$Tx_k=y_{k}$; further $\|T\|_{(E,F)}\le
2C_1(1-\delta)^{-1}.$\bull\enddemo

\proclaim{Proposition 3.6}Suppose $(E,F)$ is a pair of
K\"othe sequence spaces.  Suppose either :\newline (a)
$(E,F)$ is exponentially separated, $F$ is $r$-concave for
some $r<\infty$, and there exists $1<p<\infty$ so that $E$
has a lower $p$-estimate and $F$ has an upper
$p$-estimate.\newline or\newline (b)
$\kappa_-(E)\kappa_+(F)<1.$\newline Then $(E,F)$ has (SP) if
and only if $E$ has (RSP) and $F$ has (LSP).  \endproclaim

\demo{Proof} (a) We use Lemma 3.5 to show that $(E,F)$ has
(RSP) and $(\tilde F,\tilde E)$ has (RSP) and the result
follows.

(b) This is immediate from Lemma 3.3.\bull\enddemo

\vskip2truecm

\subheading{4.  Calder\'on couples of sequence spaces}

We now turn to calculating the K-functional for an
exponentially separated pair.

\proclaim{Lemma 4.1}Suppose $(E,F)$ is exponentially
separated.  Then there is a constant $C_2$ so that if
$\rho(a)\le t\le \rho(a+1)$ $$ K(t,x)\le
\|x_{(-\infty,a]}\|_E +t\|x_{(a,\infty)}\|_F \le
C_2K(t,x).$$ In particular, $$\|x_{(-\infty,a]}\|_E
+\rho(a)\|x_{[a,\infty)}\|_F \le C_2K(\rho(a),x).$$
Similarly if $t\le\rho(a)$ for all $a$ (in the case $\bold
J=\bold Z_+$) then $$ t\|x\|_F \le C_2K(t,x)$$ while if
$t\ge \rho(a)$ for all $a$ (when $\bold J=\bold Z_-$) then
$$ \|x\|_E \le C_2K(t,x).$$

\demo{Proof}If supp $x\subset(-\infty,a]$ then it follows
from Lemma 3.2 that $C_1K(\rho(a),x) \ge \|x\|_E.$ Similarly
if supp $x\subset[a,\infty)$ then $C_1K(\rho(a),x) \ge
\rho(a)\|x\|_F.$ Combining these statements gives the
results.\bull\enddemo

\proclaim{Theorem 4.2}Suppose $(E,F)$ is exponentially
separated and forms a Calder\'on pair.  Then $E$ satisfies
(RSP) and $F$ satisfies (LSP).  \endproclaim

\demo{Proof} First we remark that it suffices to prove the
result for $E$.  Once this is established we can apply an
order-reversal argument to get the result for $F$.  Indeed
$(\tilde F,\tilde E)$ is also exponentially separated and a
Calder\'on pair; thus $\tilde F$ has (RSP) and $F$ has
(LSP).

We will suppose that $\rho(m)\le C_02^{-n\beta }\rho(m+n)$
for $m,m+n\in\bold J$ and that $C_1$ and $C_2$ are the
constants given in Lemmas 3.2 and 4.1.

We now introduce a notion which helps in the argument.  An
admissible pair is a pair $(x,I)$ where $I$ is a finite
interval in $\bold J$ and $x$ is a positive vector with supp
$x\subset I$, max( supp $x)<\max I$, and $\|x\|_E=1.$ An
admissible family is a finite collection $\Cal
F=(x_k,I_k)_{k=1}^n$ of admissible pairs so that $(I_k)$ are
pairwise disjoint.  We define supp $\Cal F=\cup I_k.$ If
$\Cal F$ is an admissible family then we define $\Gamma(\Cal
F)$ to be the least constant $M$ so that if $(y_k)_{k=1}^n$
satisfy $\|y_k\|_E\le 1,$ supp $y_k\subset I_k$ and supp
$x_k<$ supp $y_k,$ then there exists $T\in\Cal A(E,F)$ with
$\|T\|_{(E,F)}\le M$ and $Tx_k=y_k$ for $1\le k\le n.$
Notice that since max(supp $x_k)<\max I_k$ there is ``room''
for some $y_k$ satisfying our hypotheses.  It is not
difficult to show that such a $\Gamma(\Cal F)$ is
well-defined since we can restrict the problem for each such
family to a finite-dimensional space.

We next make the remark that if $T$ is such an optimal
choice of operator then $T$ can be replaced without altering
its norm by $\sum_{k=1}^nP_{I_k}TP_{I_k}.$ Thus it can be
assumed that $Tx=0$ for any $x$ whose support is disjoint
from $\cup I_k.$ Now suppose $\Cal F$ and $\Cal G$ are two
admissible families with disjoint supports so that their
union $\Cal F\cup\Cal G$ is also admissible.  Then using the
above remark it is clear that we can simply add optimal
operators to obtain that $$ \Gamma(\Cal F\cup\Cal G)\le
\Gamma(\Cal F)+\Gamma(\Cal G).\tag 1$$

Next suppose $\Cal F$ is a single admissible pair $(x,I).$
Suppose $y$ is supported on $I$ and satisfies $\|y\|_E\le
1,$ and supp $x<$ supp $y.$ Then we can choose $x^*\in E^*$
with $\|x^*\|_{E^*}=1$ supp $x^*\subset$ supp $x$ and
$\langle x,x^*\rangle=1.$ Consider the operator $S$ defined
by $S\xi=\langle\xi,x^*\rangle y.$ Of course $\|S\|_E\le 1.$
Now suppose the maximum of supp $x$ is $a.$ Then $$ \align
\|S\xi\|_F &\le \|y\|_F \|\xi_{(-\infty,a]}\|_E \\ &\le
C_1^2\|\xi\|_F \endalign $$ where $C_1$ is the constant of
Lemma 3.2.  Hence $\Gamma(\Cal F)\le C_1^2.$ It then follows
by the addition principle (1) that if $|\Cal F|=n$ then
$\Gamma(\Cal F)\le nC_1^2.$

Now we seek to prove that $\Gamma(\Cal F)$ is bounded over
all admissible families.  Let us suppose on the contrary
that it is not.  We then can construct inductively a
sequence of admissible families $(\Cal F_n)$ for $n\in\bold
N$ and an increasing sequence of integers $(m_n)$ so that
supp $\Cal F_n\subset [-m_n,m_n]$ and $\Gamma(\Cal
F_{n+1})\ge n(n+1) + C_1^2(2m_n+2n+1).$

Now refine $\Cal F_n$ by deleting all pairs $(x,I)$ so that
$I$ intersects $[-m_n-n,m_n+n]$.  This removes at most
$2m_n+1$ pairs and creates a new admissible family $\Cal
F_n'$ so that $\Gamma(\Cal F_n')\ge n(n+1).$ The families
$\Cal F_n'$ are now disjoint.  If we write the members of
$\Cal F_n'$ in increasing order of support as
$(x_k,I_k)_{k=1}^N$ then we can define $\Cal F_{n,r}$ for
$0\le r\le n$ to be the family of all $(x_k,I_k)$ where
$k\equiv r$ mod $(n+1)$).  At least one of $\Cal F_{n,r}$
satisfies $\Gamma(\Cal F_{n,r})\ge n$ by (1).  Call this
family $\Cal G_n.$ We note that if $(x,I)$ and $(y,J)$ are
two consecutive members of $\Cal G_n$ then $I+n<J$ (since
$n$ nontrivial intervals in $\Cal F_n$ lie between $I$ and
$J.$) Furthermore there is a gap of at least $n$ between any
interval represented in $\Cal G_n$ and any interval
represented in $\Cal G_k$ for some $k<n.$

Finally let us consider the union of all $\Cal G_n$ for
$n\ge 1.$ This may be written as a sequence of admissible
pairs $(x_k,I_k)_{k\in A}$ where $A$ is one of the sets
$\bold Z,\bold Z_-,\bold Z_+$ and $I_k<I_{k+1}$ for all
$k,k+1\in A.$ Let us write $I_k=[a_k,b_k]$.  Then
$b_k<a_{k+1}$ whenever $k,k+1\in A.$ Furthermore the gaps
between the intervals tend to infinity as $|k|\to\infty.$
Precisely, if $\sigma_k=(a_{k+1}-b_k)$ then
$\lim_{|k|\to\infty}\sigma_k=\infty.$ Now let $d_k=$max
(supp $x_k$) so that $a_k\le d_k< b_k.$ Let $J_k=(d_k,b_k]$
for $k\in A.$

We now claim:

\proclaim{Claim}There exists a finite subset $A_0$ of $A$
and a constant $M$ so that if $A_1=A\setminus A_0,$ and
$(y_k)_{k\in A_1}$ is any sequence satisfying $\|y_k\|_E= 1$
and supp $y_k\subset J_k,$ then there exists $T\in\Cal
A(E,F)$ with $\|T\|_{(E,F)} \le M$ and $\|Tx_k-y_k\|_E\le
\frac12.$\endproclaim

Let us first assume the claim is established and show how
the proof is completed.  Under these hypotheses we consider
the space $\Cal Y=\ell_{\infty}(E(J_k))_{k\in A_1}$ and the
map $\Cal S:\Cal A(E,F)\to \Cal Y$ defined by $\Cal S(T) =
(P_{J_k}T_kx_k)_{k\in A_1}.$ Clearly $\|\Cal S\|\le 1$ and
it follows from the claim that if $\bold y=(y_k)_{k\in
A_1}\in \Cal Y$ there exists $T\in\Cal A(E,F)$ with
$\|T\|_{(E,F)}\le M\|\bold y\|$ and $\|\Cal S(T)-\bold
y\|\le \frac12\|\bold y\|.$ By a well-known argument from
the Open Mapping Theorem this is enough to show that $\Cal
S$ is onto and indeed if $\|\bold y\|\le 1$ then there
exists $T$ with $\|T\|_{(E,F)}\le 2M$ and $\Cal S(T)=\bold
y.$

Now suppose $\Cal G_n= \{(x_k,I_k)\}_{k\in B_n}$ where
$B_n\subset A_1.$ Then if $(y_k)_{k\in B_n}$ satisfy
$\|y_k\|_E\le 1,$ and supp $x_k<$ supp $y_k\subset I_k$ it
follows that there is an operator $T\in \Cal A(E,F)$ with
$\|T\|_{(E,F)}\le 2M$ and $P_{J_k}Tx_k=y_k.$ If we set
$T_0=\sum_{k\in B_n}P_{J_k}TP_{I_k}$ then
$\|T_0\|_{(E,F)}\le 2M$ and $T_0x_k=y_k.$ Thus $\Gamma(\Cal
G_n)\le 2M.$ Now since $A_0$ is finite we conclude that
$\Gamma(\Cal G_n)\le 2M$ for all but finitely many $n.$ This
contradicts the original construction of $\Cal G_n.$ The
contradiction shows that there is a constant $M_0$ so that
$\Gamma(\Cal F)\le M_0$ for all admissible families $\Cal
F$.  In particular if we have a finite set of finitely
supported vectors $x_1,x_2,\ldots,x_n,y_1,\ldots,y_n$ so
that supp $x_1<$ supp $y_1<\ldots<$ supp $x_n<$ supp $y_n$
and $\|x_k\|_E=1$ for all $k$ and $\|y_k\|_E\le 1$ then
there is an operator $T:E\to E$ with $\|T\|_E\le 2M_0$ and
$Tx_k=y_k.$ Hence for any $\alpha_1,\ldots,\alpha_n$ we
would have $$ \|\sum_{k=1}^n\alpha_ky_k\|_E \le
2M_0\|\sum_{k=1}^n\alpha_kx_k\| $$ and this means that $E$
has (RSP).

Thus it only remains to prove the claim.  We start by
defining a sequence $(\lambda_k)_{k\in A}$ such that
$\lambda_{k+1}-\lambda_k=\frac12\beta\sigma_k.$ We next make
some initial observations.  Let us suppose that supp
$u_k\subset I_k$ for $k\in A$ and $\|u_k\|_E= 1$ for all
$k.$ We claim that there exists a constant $C_3$ independent
of the choice of $(u_k)$ so that if $k\in A$ then

$$ \|\sum_{j\le k} 2^{\lambda_j}u_j\|_E \le
C_32^{\lambda_k}\tag 2 $$ and $$ \|\sum_{j\ge k}
2^{\lambda_j}u_j\|_F \le C_32^{\lambda_k}\|u_k\|_F.\tag 3$$

In fact (2) follows easily from the fact that if $j\le k$
then $$ \lambda_j = \lambda_k -
\frac12\beta\sum_{i=j}^{k-1}\sigma_i\le
\lambda_k-\frac12(k-j)\beta.$$

For (3) we note that if $j\ge k$, $$ \align \|u_j\|_F &\le
C_1\rho(a_j)^{-1}\\ &\le
C_1C_02^{-\beta(a_j-b_k)}\rho(b_k)^{-1}\\ &\le
C_1^2C_02^{-\beta(a_j-b_k)}\|u_k\|_F\\ &\le
C_1^2C_02^{-2(\lambda_j-\lambda_k)}\|u_k\|_F \endalign $$ so
that $$2^{\lambda_j}\|u_j\|_F \le
C_1^2C_02^{-\frac12\beta(j-k)}2^{\lambda_k}\|u_k\|_F$$ from
which (3) will follow.

In particular let us define $z=\sum_{k\in
A}2^{\lambda_k}x_k.$ The above calculations show that $z\in
E+F.$ Since $(E,F)$ is a Calder\'on couple there is a
constant $M_0=M_0(z)$ so that if $u\in E+F$ and $K(t,u)\le
K(t,z)$ for all $t$ then there exists $T\in \Cal A(E,F)$
with $\|T\|_{(E,F)}\le M_0$ and $Tz=u.$

Now suppose $(y_k)_{k\in A}$ is any sequence with
$\|y_k\|_E= 1$ and supp $y_k\subset J_k.$ We set $v=\sum
2^{\lambda_k}y_k \in E+F.$ We turn to comparing $K(t,v)$
with $K(t,z).$ Let us note first that for every $k\in A$ we
have $$ \|y_k\|_F \le C_1\rho(d_k)^{-1}\le C_1^2\|x_k\|_F.$$

If $t$ satisfies $t\le \rho(a_k)$ for all $k$ then we must
have $A=\bold Z_+$ and we make the estimate $$ K(t,v)\le
t\|v\|_F \le C_3t2^{\lambda_0}\|y_0\|_F \le
C_3C_1^2t2^{\lambda_0}\|x_0\|_F$$ so that by Lemma 3.3, $$
K(t,v) \le C_3C_2C_1^2K(t,z).$$

Similarly if $t\ge \rho(b_k)$ for all $k\in A$ then we can
have $A=\bold Z_-$ and we make a similar estimate $$ K(t,v)
\le \|v\|_E \le C_32^{\lambda_{-1}}\le
C_3\|z_{[a_{-1},b_{-1}]}\|_E \le C_2C_3K(t,z).$$

In the other cases we first consider the case when
$\rho(n)\le t\le \rho(n+1)$ for some $n$ in an interval
$[a_k,d_k).$ Then $K(t,x_k) \ge tC_1^{-1}\rho(d_k)^{-1}$ by
Lemma 3.2.

Hence $$ \align K(t,v_{[a_k,\infty)}) &\le
C_3t\|v_{[a_k,a_{k+1}]}\|_F\\ &=
C_3t2^{\lambda_k}\|y_k\|_F\\ &\le C_3C_1
t2^{\lambda_k}\rho(d_k)^{-1}\\ &\le C_3C_1^2
2^{\lambda_k}K(t,x_k)\\ &\le C_3C_1^2 K(t,z).  \endalign $$
If $k$ is the initial element of $A$ we are done.
Otherwise:  $$ \align K(t,v_{(-\infty,a_k)}) &\le C_3
2^{\lambda_{k-1}}\\ &\le C_3\|z_{(-\infty,a_k)}\|_E \\ &\le
C_3 K(t,z).  \endalign $$ Combining in this case we have
$K(t,v)\le CK(t,z)$ for some constant $C$ depending only on
$C_1,C_2$ and $C_3.$

For the final case, we can suppose there exists $n$ not in
any interval $[a_k,d_k)$ and such that $\rho(n)\le t\le
\rho(n+1);$ it may also be assumed that there exists $k\in
A$ with $k+1\in A$ and $d_k\le n\le a_{k+1}.$ Then by Lemma
3.3, $$ \|z_{(-\infty,n]}\|_E + t\|z_{(n,\infty)}\|_F \le
C_2K(t,z).$$ Now $$ \|v_{(-\infty,b_k]}\|_E \le C_3
2^{\lambda_k}\le C_3\|z_{(-\infty,n]}\|_E.$$ Also $$
\|v_{(b_k,\infty)}\|_F \le C_3
2^{\lambda_{k+1}}\|y_{k+1}\|_F \le
C_1^2C_3\|z_{(n,\infty)}\|_F.$$

Thus combining all the cases there exists $C_4$ independent
of $(y_k)$ so that $K(t,v)\le C_4K(t,z)$.  Hence there is an
operator $T\in\Cal A(E,F)$ with $Tz=v$ and $\|T\|_{(E,F)}\le
C_4M_0.$

Now for fixed $k\in A$ assume first that $k$ is not the
initial element of $A$.  Then $$ \|(z_{(-\infty,a_k)}\|_E
\le C_32^{\lambda_{k-1}}\le
C_32^{-\frac12\beta\sigma_{k-1}}2^{\lambda_k}.$$ Thus we
have that $$ \|T(z_{(-\infty,a_k)})\|_E \le
C_4C_3M_02^{-\frac12\beta\sigma_{k-1}}2^{\lambda_k}.$$

If $k$ is not the final element, $$ \align \|
z_{(b_k,\infty)}\|_F &\le C_3 2^{\lambda_{k+1}}\|y_{k+1}\|_F
\\ &\le C_1C_3 2^{\lambda_{k+1}}\rho(a_{k+1})^{-1}\\ &\le
C_0C_1C_3 2^{\lambda_{k+1}-\beta\sigma_k}\rho(b_k)^{-1}.
\endalign $$ Thus if $f=T(z_{(b_k,\infty)})_{(d_k,b_k]}$
then $$ \align \|f\|_E &\le C_1\rho(b_k)\|f\|_F\\ &\le
C_1C_4M_0\rho(b_k) \|z_{(b_k,\infty)}\|_F\\ &\le
C_0C_1^2C_3C_4M_02^{-\frac12\beta\sigma_k}2^{\lambda_k}.
\endalign $$

It follows that if $k$ is not an initial or final element of
$A,$ $$ \|y_k- P_{J_k}Tx_k\|_E \le
C_5M_02^{-\frac12\beta\tau_k}$$ where $C_5$ is a constant
depending only on $E$ and $F$ and
$\tau_k=\min(\sigma_{k-1},\sigma_k).$ Now if we set
$S=\sum_{k\in A}P_{J_k}TP_{I_k}$ then $\|S\|_{(E,F)}\le
C_4M_0$.  Further if we let $A_1$ be the set of $k\in A$ so
that $k$ is not an initial or final element and
$C_5M_0^22^{-\frac12\beta\tau_k}<1/2$ then $A_0=A\setminus
A_1$ is finite and $ \|Sx_k-y_k\| < 1/2$ for $k\in A_1.$
Thus the claim is established and the proof is
complete.\bull\enddemo

\proclaim{Lemma 4.3}Suppose $(E,F)$ satisfies (RSP).  Then
there is a constant $C$ so that if $0\le x,y\in E+F$ and
$\|y_{(-\infty,a]}\|_E \le \|x_{(-\infty,a]}\|_E$ for all
$a\in\bold J$ then there exists a positive $T\in\Cal A(E,F)$
with $\|T\|_{(E,F)}\le C$ and $Tx=y.$\endproclaim

\demo{Proof}By applying the argument of Lemma 2.6 we deduce
from (RSP) the existence of a constant $C_0$ so that if $A$
is an interval in $\bold Z$, $(J_k)_{k\in A}$ is a
collection of finite intervals in $\bold J$ with
$J_k<J_{k+1}$ whenever $k,k+1\in A$ and $(x_k)_{k\in
A},(y_k)_{k\in A}$ are positive and satisfy supp $x_k,$ supp
$y_k\subset J_k$ and $\|y_{k+1}\|_E\le \|x_k\|_E$ for
$k,k+1\in A$ then there is a positive matrix operator $T$
with $\|T\|_{(E,F)}\le C_0$ and $Tx_k=y_{k+1}$ whenever
$k,k+1\in A.$

Let us prove the lemma when $x,y$ have disjoint supports.
We first define a function $\sigma:\bold Z\to\bold
Z\cup\{\pm\infty\}$ by setting $\sigma(k)=-\infty$ if
$k<\bold J$, $\sigma(k)=\infty$ if $k>\bold J$ and otherwise
$\sigma(k)$ is the greatest $j\in \bold Z\cup\{-\infty\}$ so
that $\|x_{(-\infty,k]}\|_E \ge 4^j.$

Let $I_0=\{k\in\bold J:\sigma(k)>\sigma(k-1)\}.$ We then let
$I$ be the subset of $I_0$ of all $k$ so that for every
$n\in I_0$ with $n<k$ then $\|x_{(-\infty,n]}\|_E \le
\frac12 \|x_{(-\infty,k]}\|_E.$

We can now index $I$ as $(a_n)_{n\in A}$ where $A$ is an
interval in $\bold Z$ which can be assumed to have $0$ as
its initial element if $I$ is bounded below.

We now define $B$ to be $\bold Z_-$ when $\inf_{k\in\bold
Z}\sigma(k)>-\infty$ and to be empty otherwise.  We will
only need to introduce $B$ in the case when $\lim_{a\to
-\infty}\|x_{(-\infty,a]}\|_E >0.$ If $B$ is nonempty then
$I$ is bounded below and there exists a greatest $\lambda$
so that $\|x_{(-\infty,k]}\|_E \ge 4^{\lambda}$ for every
$k$ (in this case $\bold J$ cannot be bounded below).  We
must have $\|x_{(-\infty,k]}\|_E <4^{\lambda+1}$ whenever
$k<a_0.$ It follows that we may pick $a_{-1}$ so that
$\|x_{(a_{-1},a_0]}\|_E \ge 4^{\sigma(a_0)-1}$ and then
inductively $a_{-n}$ so that $\|x_{(a_{-n},a_{-(n-1)}]}\|_E
\ge 4^{\lambda-1}.$

In this way we define $(a_n)_{n\in A\cup B}.$ We now let
$x_n=x_{(a_{n-1},a_n]}$ and $y_n=y_{(a_{n-1}a_n]}$ if $n\in
A\cup B$ is not the initial element of $A\cup B$; if $n=0$
is the initial element we let $x_0=x_{(-\infty,a_0]}$ and
$y_0=y_{(-\infty,a_0]}.$ If $n$ is the final element of
$A\cup B$ we set $y_{n+1}=y_{(a_n,\infty]}.$ We may now
verify that $\sum_{n\in A\cup B}x_n\le x.$ We also claim
that $\sum_{n\in A\cup B}y_{n+1}=y.$ If $A\cup B=\bold Z$
this is clear.  If $A\cup B=(-\infty,n]$ for some $n$ it
follows from our definition of $y_{n+1}.$ If $A\cup B$ is
bounded below (by 0) then $B$ is empty and hence
$\sigma(a_0-1)=-\infty.$ Thus $y_0=0$ and we obtain our
claim easily.

We first prove that if $n,n+1\in A$ then
$\sigma(a_{n+1}-1)\le \sigma(a_n)+1.$ If not there exists a
first $k_1$ so that $\sigma(k_1)\ge \sigma(a_n)+1$ and a
first $k_2$ so that $\sigma(k_2)\ge \sigma(a_n)+2$ and
$a_n<k_1\le k_2\le a_{n+1}-1.$ Then $k_1,k_2$ are in
$I_0\setminus I.$ Thus $\|x_{(-\infty,a_n]}\|_E >\frac12
\|x_{(-\infty,k_1]}\|_E.$ The equality $k_1=k_2$ would
entail $\|x_{(-\infty,k_1]}\|_E \ge 4^{\sigma(a_n)+2}$ and
thus $\|x_{(-\infty,a_n]}\|_E > 4^{\sigma(a_n)+1}$ which
contradicts the definition of $\sigma(a_n).$ Thus $k_1<k_2$
and we conclude also that $\|x_{(-\infty,k_1]}\|_E >\frac12
\|x_{(-\infty,k_2]}\|_E$ so that $\|x_{(-\infty,a_n]}\|_E
>\frac14 \|x_{(-\infty,k_2]}\|_E$ which implies the absurd
conclusion $\sigma(k_2)\le \sigma(a_n)+1.$ Thus, as claimed,
$\sigma(a_{n+1}-1) \le \sigma(a_n)+1.$

The same argument shows that if $A$ is bounded above then if
$k>a_n$ we must have $\sigma(k)\le \sigma(a_n)+1.$

Now if $n,n+1\in A$ we can argue that since $x,y$ have
disjoint supports that $y_{n+1}$ is supported on
$(a_n,a_{n+1})$ and thus $\|y_{n+1}\|_E \le
\|x_{(-\infty,a_{n+1}-1]}\|_E \le 4^{\sigma(a_n)+2}.$
Similarly, let $n$ be the last element of $A$.  Then for all
$k>a_n$ $\|y_{(-\infty,k]}\|_E \le \|x_{(-\infty,k]}\|_E\le
4^{\sigma(k)+1}\le 4^{\sigma(a_n)+2}.$ Thus $\|y_{n+1}\|_E
\le 4^{\sigma(a_n)+2}.$

On the other hand, if $n$ is not the initial element of $A$,
$$\|x_n\|_E \ge
\|x_{(-\infty,a_n]}\|_E-\|x_{(-\infty,a_{n-1}]}\|_E \ge
\frac12 4^{\sigma(a_n)}.$$ If $n=0$ is the initial element,
we either have, if $B=\emptyset$, $x_0=x_{(-\infty,a_0]}$ so
that $\|x_0\|\ge 4^{\sigma(a_0)}$ or if $B\neq \emptyset$
then $\|x_1\| \ge 4^{\sigma(a_0)-1}.$ In all such cases, if
$n\in A$ we have $\|y_{n+1}\|_E\le 4^3 \|x_n\|_E.$

Next suppose $n,n+1\in B.$ Then $\|y_{n+1}\|_E \le
\|x_{(-\infty,a_{n+1})}\|_E \le 4^{\lambda+1}$ while
$\|x_n\| \ge 4^{\lambda-1}.$ Thus $\|y_{n+1}\|_E\le
4^2\|x_n\|_E.$

Finally, consider the case $n=-1\in B$ and $n+1=0\in A$.
Then since $a_0$ is in the support of $x$ we have
$\|y_{n+1}\|_E \le \|x_{(-\infty,a_0-1)}\|_E <
4^{\lambda+1}.$ However $\|x_n\|_E \ge 4^{\lambda-1}$ so
that $\|y_{n+1}\|_E\le 4^2\|x_n\|_E.$

Combining all cases, we conclude that there is a positive
operator $T$ with $\|T\|_{(E,F)}\le 4^3C_0$ so that
$Tx_n=y_{n+1}$.  Now it is clear that $\sum_{n\in A\cup
B}x_n\le x$ while $\sum_{n\in A\cup B}y_{n+1}=y.$ Thus if
$S=\sum_{n\in A\cup B}P_{\text{supp }y_{n+1}}TP_{\text{supp
}x_n}$ then $\|S\|_{(E,F)}\le 64C_0$ and $Sx=y.$ Thus the
lemma is established in the case when $x$ and $y$ have
disjoint supports.

For the general case we let $I=\{n\in\bold J:y_n >2x_n\}.$
Let $J=\bold J\setminus I.$ Then set $u=x_J$ and $v=y_I.$
For any $k\in \bold J$ we have $\|x_{I\cap(-\infty,k]}\|_E
\le \frac12\|y_{(-\infty,k]}\|_E \le
\frac12\|x_{(-\infty,k]}\|_E.$ Thus $\|u_{(-\infty,k]}\|_E
\ge \frac12\|x_{(-\infty,k]}\|_E.$ Hence there is a positive
operator $S$ in $\Cal A(E,F)$ with $Su=v$ and
$\|S\|_{(E,F)}\le 128C_0.$ On the other hand $y_J\le 2x$ and
so there is a multiplication operator $V\in\Cal (E,F)$ with
$\|V\|_{(E,F)}\le 2$ and $Vx=y_J.$ Finally the operator
$T=SP_J + V$ establishes the lemma.\bull\enddemo

\proclaim{Lemma 4.4}Suppose $(E,F)$ is exponentially
separated and satisfies (SP).  Then there exists a constant
$C$ so that if $0\le x,y\in E+F$ and $K(t,y)\le K(t,x)$ for
all $t\ge 0$ then there exists a positive matrix $T\in\Cal
A(E,F)$ with $\|T\|_{(E,F)}\le C$ and $Tx=y.$

\demo{Proof}It follows from Lemma 4.1 that there is a
constant $C_2$ so that for all $a\in\bold J$ we have,
whenever $K(t,y)\le K(t,x)$ for all $t\ge 0$, $$
\max(\|y_{(-\infty,a]}\|_E,\rho(a)\|y_{[a,\infty)}\|_F)\le
2C_2\max(
\|x_{(-\infty,a]}\|_E,\rho(a)\|x_{(a,\infty)}\|_F).$$ Thus
for every $a$ either $$ \|y_{-\infty,a]}\|_E \le
2C_2\|x_{(-\infty,a]}\|_E \tag 4$$ or $$
\|y_{[a,\infty)}\|_F \le 2C_2\|x_{[a,\infty)}\|_F.  \tag 5$$
Let $J_1$ be the set of $a$ so that (4) holds and let
$J_2=\bold J\setminus J_1.$ Since $(E,F)$ has (RSP) we can
apply Lemma 4.3 to deduce the existence of a positive matrix
$T_1$ with $\|T_1 \|_{(E,F)} \le C_3$ where $C_3$ depends
only $(E,F)$ and $T_1x=y_{J_1}.$ Similarly since $(E,F)$ has
(LSP) we can find a positive matrix $T_2$ with
$\|T_2\|_{(E,F)}\le C_3$ so that $T_2x=y_{J_2}.$ Then
$(T_1+T_2)(x)=y.$ \bull\enddemo

\proclaim{Theorem 4.5}Let $(E,F)$ be a pair of K\"othe
sequence spaces.  Suppose either:\newline (a)
$\kappa_-(E)\kappa_+(F)<1.$ \newline or\newline (b) $(E,F)$
is exponentially separated, $F$ is r-concave for some
$r<\infty$ and there exists $p$ with $1\le p<\infty$ so that
$E$ has a lower $p$-estimate and $F$ has an upper
$p$-estimate.\newline Then $(E,F)$ is a (uniform) Calder\'on
couple if and only if $E$ has (RSP) and $F$ has
(LSP).\endproclaim

\demo{Proof}This is an immediate deduction from Proposition
3.6, Theorem 4.2 and Lemma 4.4.\bull\enddemo

\demo{Remark}Note that in fact Lemma 4.4 implies that under
these circumstances if $K(t,y)\le K(t,x)$ for all $t$ and
$x,y\ge 0$ then there is a positive operator $T$ with
$\|T\|_{(E,F)}<\infty$ and $Tx=y.$\enddemo

The following theorem is similar to results of Cwikel and
Nilsson [18].

\proclaim{Theorem 4.6} Let $E,F$ be symmetric sequence
spaces on $\bold Z_+$ and suppose $(E,F(w))$ is a Calder\'on
pair for a weight sequence $w=(w_n)$.  Then either $F(w)=F$
(i.e.  $0<\inf w_n\le \sup w_n<\infty$) or $E=\ell_p,\
F=\ell_q$ for some $1\le p,q\le\infty.$\endproclaim

\demo{Proof}If $(w_n)$ is unbounded we can pass to a
subsequence satisfying $w_{n_k}>2w_{n_{k-1}}$.  Then the
pair $(E,F(w_{n_k}))$ is a Calder\'on pair and we can apply
Theorem 4.2 to get that $F$ has (RSP) and $E$ has (LSP).  An
application of Proposition 2.3 gives the result.  If
$(w_n^{-1})$ is unbounded we can argue
similarly.\bull\enddemo

\subheading{5.  Calder\'on couples of r.i. spaces}

Let $\Omega$ denote one of the sets $[0,\infty)$, $[0,1]$
and $\bold N$.  Let $\bold J$ be the set $\bold Z$, $\bold
Z_-,$ or $\bold Z_+$ respectively.  If $X$ is an r.i. space
on $\Omega$ (or a symmetric sequence space if $\Omega=\bold
N)$ we will associate to $X$ a K\"othe sequence space $E_X$
on $\bold J.$ To do this let $e_n,\ n\in\bold J$ be defined
by $e_n=\chi_{[2^n,2^{(n+1)})}.$ We then define for
$x\in\omega(\bold J),$ $$ \|x\|_{E_X} = \|\sum_{k\in\bold
J}x(k)e_k\|_X.$$ (Here we use $e_k$ with a dual meaning as
both the canonical basis element of $\omega(\bold J)$ and as
an element of $X(\Omega).$) We observe that $E_X$ regarded
as a subspace of $X$ is 1-complemented by the natural
averaging operator.  Notice also that $E_{X^*}= E_X^*(2^n)$
is a weighted version of $E_X^*$.  We also note that on
$E_X$ we can compute $\|\tau_n\|_{E_X}\le \|D_{2^n}\|_X$
where $D_s$ is the natural dilation operator.  Furthermore
it is easy to see that for $f\in X$ we have $D_{2^n}f^*\le
\tau_{n+1}Pf^*$ where $P$ is the natural averaging
projection of $X$ onto $E_X$; thus $\|D_{2^n}\|_X \le
\|\tau_{n+1}\|_{E_X}.$ Thus $\kappa_+(E_X)= 2^{1/p_X}$ and
$\kappa_-(E_X) =2^{-1/q_X}$ where $p_X$ and $q_X$ are the
Boyd indices of $X$.

We now show how to build examples of r.i. spaces from
sequence spaces.  To keep the notation straight we prove our
results for the case of function spaces $\Omega=[0,1]$ or
$\Omega=[0,\infty)$.  However simple modifications give the
analogous results for sequence spaces.

\proclaim{Proposition 5.1}Let $E$ be a K\"othe sequence
space on $\bold J.$ Then:\newline (1) If $\kappa_+(E)<2$
there is an r.i. space $X=X(\Omega)$ so that $\|f\|_X$ is
equivalent to $\|\sum_{n\in\bold J}f^*(2^n)e_n\|_E$.\newline
(2) If $\kappa_-(E)<1 \le\kappa_+(E)<2,$ and $X$ is an r.i.
space so that $\|f\|_X$ is equivalent to $\|\sum_{n\in\bold
J}f^*(2^n)e_n\|_E$ then $E_X=E$ (up to equivalence of norm).
\endproclaim

\demo{Proof} (1) We define $X$ to be the set of measurable
functions on $\Omega$ such that $$ \|f\|_X =
\|\sum_{n\in\bold J}f^*(2^n)e_n\|_{E}<\infty.$$ We show that
the functional $\|f\|_X$ is equivalent to a norm by
computing $\|f^{**}\|_X$ where
$f^{**}(t)=\frac1t\int_0^tf^*(s)ds.$ Then $$ f^{**}(2^n) \le
2^{-n}\sum_{k<n}2^kf^*(2^k).$$ Thus $$ \|f^{**}\|_X \le
\sum_{j=1}^{\infty}2^{-j}\|\tau_{j}(\sum_{n\in\bold
J}f^*(2^n)e_n) \|_E.$$ Thus since $\kappa_+(E)<2$ there is a
constant $C_1$ so that $\|f^{**}\|_X \le C_1\|f\|_X.$ Since
$f\to \|f^{**}\|_X$ is plainly an r.i. norm and the set
$\{f:\|f^{**}\|_X\le 1\}$ is closed in measure it is clear
that $X$ is an r.i. space.

(2) Let $\|\ \|_X$ denote the quasinorm induced by $E.$ We
remark that it follows from (1) that there exists a constant
$C_2$ so that for $f\in X$ we have $\|f^{**}\|_X\le
C_2\|f\|_X$.  Now, considering the $E_X$-quasinorm induced
on $\omega(\bold J)$ it is clear that if $x$ is a
nonincreasing sequence then $\|x\|_{E_X}=\|x\|_E.$ In
general we note that if $f\in E_X$ then for some
$C_3=\sum_{j\ge 0}\|\tau_{-j}\|_E$ since $\kappa_-(E)<1,$ $$
\|\max_{j\ge 0}|\tau_{-j}f|\|_E \le C_3\|f\|_E$$ so that $$
\|f\|_{E_X} \le C_1\|f\|_E.$$ For the converse direction we
observe that if $f\in X$ it is trivial that $\|D_{2^j}f\|_X
\le \|\tau_j\|_E\|f\|_X.$ Then $$ \align \|f\|_E &\le
\|\max_{j\ge 0}|\tau_{-j}f|\|_E\\ &= \|\max_{j\ge
0}|\tau_{-j}f|\|_X\\ &\le \| \sum_{j\ge 0} D_{2^{-j}}|f|\|_X
\\ &\le \sum_{j\ge 0} \|D_{2^{-j}}f^{**}\|_X \\ &\le C_3
\|f^{**}\|_X\\ &\le C_2C_3\|f\|_X.  \endalign $$ Thus $E_X$
is (up to equivalence of norm) identical with $E.$
\bull\enddemo

\demo{Remark} It follows from the above Proposition that
there is a natural one-one correspondence between r.i.
spaces $X$ with Boyd indices satisfying $1<p_X\le
q_X<\infty$ and sequence spaces $E$ on $\bold J$ with
$\kappa_-(E)<1\le \kappa_+(E)<2$ determined by $E=E_X.$
Under this correspondence if $1<p<\infty$ an r.i. space $X$
with $q_X<\infty$ is a Lorentz space (of order $p$) if and
only if $E_X$ is a weighted $\ell_p-$ space.  For if
$$\|f\|_X =
(\int_0^{\infty}(f^*(t)w(t))^p\frac{dt}{t})^{1/p}$$ where
$w$ is an increasing function satisfying $1<\inf w(2t)/w(t)
\le \sup w(2t)/w(t) <\infty$ then the above Proposition
shows that $E_X= \ell_p(w_n),$ where $w_n=w(2^n).$
Conversely if $E_X$ is an $\ell_p-$space then
$E_X=\ell_p(w_n)$ where the assumption that $q_X<\infty$
enables us to assume $\inf w_{n+1}/w_n>1.$ If we define
$w(t)=w_n$ whenever $2^{n-1}<t\le 2^n$ then it is easy to
see that $X$ is a Lorentz space.

We now prove the elementary:

\proclaim{Proposition 5.2}Let $(X,Y)$ be a pair of r.i.
spaces on $\Omega$.  Then $(X,Y)$ is a Calder\'on couple if
and only if $(E_X,E_Y)$ is a Calder\'on couple.\endproclaim

\demo{Proof}By using the averaging projection it is clear
that if $(X,Y)$ is a Calder\'on couple then so is
$(E_X,E_Y).$ In fact it is trivial to see that for $f\in
E_X+E_Y$ we have $K(t,f;X,Y)=K(t,f;E_X,E_Y).$ Thus if
$K(t,g;E_X,E_Y)\le K(t,f;E_X,E_Y)$ for all $t$ there exists
$T\in \Cal A(X,Y)$ so that $Tf=g$.  If $P$ is the averaging
projection then $PT\in\Cal A(E_X,E_Y)$ and $PTf=g.$

Conversely suppose $(E_X,E_Y)$ is a Calder\'on couple.
Suppose $f,g\in X+Y$ and $K(t,g;X,Y)\le K(t,f;X,Y)$ for all
$t\ge 0.$ We then observe that if $G=\sum_{n\in\bold
J}g^*(2^n)e_n$ and $F=\sum_{n\in\bold J}f^*(2^n)e_n$ then
$g^*\le G\le D_2g^*$ and $f^*\le F\le D_2f^*.$ and $$
K(t,G;X,Y) \le K(t,D_2g^*;X,Y) \le 2K(t,f;X,Y)\le
2K(t,F;X,Y).$$ Since $F,G$ are in $E_X+E_Y$ we can deduce
the existence of $T\in\Cal A(E_X,E_Y)$ with $TF=G.$ Now
since $F\le D_2f^*$ and $g^*\le G$ it is clear that there
exists $S\in \Cal A(X,Y)$ with $Sf=g$.\bull\enddemo

\demo{Remarks}It now follows that every pair of Lorentz
spaces whose Boyd indices are finite is a Calder\'on couple,
since every pair of weighted $\ell_p-$spaces is a Calder\'on
couple (cf.  [36], [13]); this result is due to Cwikel [14]
and Merucci [30] for certain special cases.

We introduce the following definitions.  We say $X$ is {\it
stretchable} if $E_X$ has (RSP) and we say that $X$ is {\it
compressible} if $E_X$ has (LSP).  If $X$ both stretchable
and compressible, we say that $X$ is {\it elastic.} It is
immediate from Proposition 2.1 that $X$ is stretchable if
and only if $X^*$ is compressible and {\it vice versa}; thus
elasticity is a self-dual property.  We remark that we have
no example of a stretchable (or compressible) space which is
not already elastic.  In fact we shall see that for Orlicz
spaces these concepts do indeed coincide.

\proclaim{Theorem 5.3} Let $(X,Y)$ be a pair of r.i. spaces
on $\Omega$ whose Boyd indices satisfy $p_Y>q_X.$ Then
$(X,Y)$ is a Calder\'on couple if and only if $X$ is
stretchable and $Y$ is compressible.  \endproclaim

\demo{Proof}As $\kappa_-(E_X)= 2^{-1/q_X}$ and
$\kappa_+(E_Y)= 2^{1/p_Y}$ we have
$\kappa_-(E_X)\kappa_+(E_Y)<1$ and so the theorem is
immediate from Theorems 4.5 and 5.2.\bull\enddemo

If one space is $L_{\infty}$ we can do rather better.

\proclaim{Theorem 5.4} Let $X$ be an r.i. space on
$\Omega=[0,1]$ or $\Omega=[0,\infty).$ Then $(X,L_{\infty})$
is a (uniform) Calder\'on couple if and only if $X$ is
stretchable.  Similarly if $X$ is a symmetric sequence space
then $(\ell_{\infty},X)$ is a (uniform) Calder\'on couple if
and only if $X$ is stretchable.  \endproclaim

Before proving Theorem 5.4 we state a result which has a
very similar proof.  We remark that Theorem 5.5 only
improves on Theorem 5.3 under the assumption that
$p_Y=p=q_X$ since the case $p_Y<q_X$ is already covered.

\proclaim{Theorem 5.5}Suppose $(X,Y)$ is a couple of r.i.
spaces on $\Omega$ so that for some $1\le p<\infty$ $X$ is
$p$-concave and $Y$ is $p$-convex and suppose also that $Y$
is $r$-concave for some $r<\infty.$ Then $(X,Y)$ is a
(uniform) Calder\'on couple if and only if $X$ is
stretchable and $Y$ is compressible.\endproclaim

\demo{Proofs of Theorems 5.4 and 5.5} Theorem 5.4
corresponds to the case $p=\infty,$ and $Y=L_{\infty}.$ We
can and do assume that the $p$-convexity constant of $Y$ and
the $p$-concavity constant of $X$ are both equal to one.
Under this hypothesis it is easy to see that, when
$p<\infty$ $2^{-k/p}\|e_k\|_X$ is increasing and
$2^{-k/p}\|e_k\|_Y$ is decreasing.  Thus for $p\le \infty,$
$\rho(k)=\|e_k\|_X/\|e_k\|_Y$ is an increasing function and
$\rho(k+1)\le 2\rho(k)$ whenever $k,k+1\in\bold J.$ Then for
$k\in \bold J$ we let $I_k=\{n\in\bold J:2^k< \rho(n)\le
2^{k+1}\}.$

Before continuing let us make remark which we use several
times in the proof.  Assuming $p<\infty$ suppose $f,g$ are
two finitely supported functions in $E_X$ which satisfy
$\|f\|_p=\|g\|_p$ and $$ \int_0^t (f^*(s))^pds \le \int_0^t
(g^*(s))^pds$$ for every $t\ge 0.$ Then we have the
inequalities $\|f\|_Y \le \|g\|_Y$ and $\|f\|_X\le \|g\|_X.$
In fact it follows from a well-known lemma of Hardy,
Littlewood and Polya, [19], [25], p.124, that $|f^*|^p$ is
in the convex hull of the set of all rearrangements of
$|g^*|^p$; this can be proved by partitioning the supports
of $f^*,g^*$ into finitely many sets of equal measure.  The
assertion is then a direct consequence of the definitions of
$p$-convexity and $p$-concavity.

We make some initial remarks which will be needed in both
directions of the proof.  Each set $I_k$ is an interval
(possibly infinite) or is empty.  The set of $k$ so that
$I_k$ is nonempty is an interval $A$.  Let $E(I_k)$ be the
linear span of $(e_n:n\in I_k)$ when $k\in A.$ We state the
following Lemma.

\proclaim{Lemma 5.6} If $f,g\in E(I_k)$ then, under the
hypotheses of Theorem 5.5, $$ \align \|f\|_X\|g\|_p &\le
2\|f\|_p\|g\|_X\\ \|f\|_Y\|g\|_p &\le 2\|f\|_p\|g\|_X
\endalign $$ where $\|\ \|_p$ denotes the usual $L_p$-norm,
so that $\|\sum \alpha_ke_k\|_p =(\sum 2^k|\alpha_k|^p
)^{1/p}$.\newline Under the hypotheses of Theorem 5.4, we
have $$ \|f\|_X\|g\|_{\infty} \le 4\|f\|_{\infty}\|g\|_X.$$
\endproclaim

\demo{Proof} In fact suppose $f,g\in [e_n:a\le n\le b]$
where $a,b\in I_k,$ and that neither is zero.  We may
observe that for all $t\ge 0$ we have
$$2^{-a}\int_0^t(e_a^*(s))^p)ds \ge \|f\|_p^{-p}\int_0^t
(f^*(s))^pds \ge 2^{-b}\int_0^t (e_b^*(s))^pds$$ with
similar inequalities for $g$.  It thus follows from the
remarks above that $$ 2^{-a/p}\|e_a\|_Y \ge
\|f\|_p^{-1}\|f\|_Y \ge 2^{-b/p}\|e_b\|_Y.$$ Similarly $$
2^{-a/p}\|e_a\|_X \le \|f\|_p^{-1}\|f\|_X \le
2^{-b/p}\|e_b\|_X.$$ There are similar inequalities for $g.$
Since $2^k<\rho(a)\le\rho(b)\le 2^{k+1}$ $$
2^{-b/p}\|e_b\|_X \le 2^{k+1-b/p}\|e_b\|_Y \le
2^{k+1-a/p}\|e_a\|_Y \le 2.2^{-a/p}\|e_a\|_X.$$ Combining
these we see that $$ \|f\|_p^{-1}\|f\|_X \le
2\|g\|_p^{-1}\|g\|_X$$ and $$ \|f\|_p^{-1}\|f\|_Y \le
2\|g\|_p^{-1}\|g\|_Y$$ whence the claimed inequalities
follow.  For the last part, we observe that
$$\|f\|_{\infty}\|e_a\|_X \le \|f\|_X \le
\|f\|_{\infty}\|e_a+\cdots+e_b\|_X \le
2\|f\|_{\infty}\|e_b\|_X$$ and proceed
similarly.\bull\enddemo

We draw immediately the conclusion that if $A$ is finite (so
that $\rho$ is bounded) then both $X$ and $Y$ coincide with
$L_p(\mu)$ and there is nothing to prove.  In other cases at
most one $I_k$ is infinite.  We write $I_k=[a_k,b_k]$ if
$I_k$ is finite and $I_k=[a_k,\infty)$ or
$I_k=(-\infty,b_k]$ if $I_k$ is infinite.  Let $A_0$ be the
set of $k$ so that $k-1$ and $k+1\in A$.  We define a set
$J$ by taking one point $d_k$ from each $I_k$ for $k\in A$.
We introduce the sequence spaces $F_X$ and $F_Y$ modelled on
$A$ by setting $\|x\|_{F_X} =\|\sum_{k\in A}
2^{-d_k/p}x(k)e_{d_k}\|_X$ and $\|x\|_{F_Y} =\|\sum_{k\in
A}2^{-d_k/p}x(k)e_{d_k}\|_Y.$ In the case $p=\infty$ we
define $\|x\|_{F_X}= \|\sum_{k\in A}x(k)e_{d_k}\|_X.$

\proclaim{Lemma 5.7}Under the hypotheses of Theorem 5.5,
suppose $E_Y(J)$ has (LSP).  Then there is a constant $C_0$
so that if $f\in E_Y$ then $\|f\|_Y$ is $C_0-$equivalent to
$\| (\|f_{I_k}\|_p)\|_{F_Y}.$ \endproclaim

\demo{Proof}It suffices to prove such an equivalence if
$f\in E_Y$ satisfies $f_{I_k}=0$ for $k\notin A_0,$ since
there are most two values of $k\notin A_0$ and Lemma 5.6
shows that the $Y$-norm on each such $E(I_k)$ is equivalent
to the $L_p-$norm.  Next observe that for such $f$ if
$g=\sum_{k\in A_0}2^{-d_{k+1}/p}\|f_{I_k}\|_pe_{d_{k+1}}$
then for all $t\ge 0,$ $\int_0^t(g^*(s))^pds \le
\int_0^t(f^*(s))^pds.$ Thus we have immediately by the
$p$-convexity and rearrangement-invariance of $Y,$
$\|g\|_Y\le \|f\|_Y.$ Similarly if $h=\sum_{k\in
A_0}2^{-d_{k-1}/p}\|f_{I_k}\|_pe_{d_{k-1}}$ then $\|h\|_Y\ge
\|f\|_Y.$ Next let $\tilde f=\sum_{k\in
A_0}2^{-d_k/p}\|f_{I_k}\|_pe_{d_k}.$ We complete the proof
by showing that for some $C,$ $\|h\|_Y\le C\|\tilde f\|_Y$
and $\|\tilde f\|_Y\le C \|g\|_Y.$ Once this is done it will
be clear that $\|f\|_Y$ is actually equivalent to $\|\tilde
f\|_Y$ as claimed.

The proofs of these statements are essentially the same, so
we concentrate on the first.  Note that $$ \align
2^{-d_{k-1}/p}\|e_{d_{k-1}}\|_Y &\le
2^{-(k-1)-d_{k-1}/p}\|e_{d_{k-1}}\|_X\\ &\le
2^{-(k-1)-d_k/p}\|e_{d_k}\|_X \\ &\le
2^{1-d_k/p}\|e_{d_k}\|_Y \endalign $$ and so if $C$ is the
(LSP) constant of $E_Y(J)$ we have $\|h\|_Y \le 2C\|\tilde
f\|_Y.$ Similarly $\|\tilde f\|_Y \le
2C\|g\|_Y.$\bull\enddemo

In a very similar way, exploiting the $p$-concavity of $X$
one has,

\proclaim{Lemma 5.8} Suppose $E_X(J)$ has (RSP).  Then there
is a constant (which we also name $C_0$) so that if $f\in
E_X$ then $\|f\|_X$ is $C_0-$equivalent to $\|
(\|f_{I_k}\|_p)\|_{F_X}.$ \endproclaim

\demo{Sketch}First consider the case of Theorem 5.5.  We
assume $f\in E_X$ is finitely supported.  Proceed as in
Lemma 5.7, defining $g,h,\tilde f$ as before.  In this case
we have that $\|g\|_X\ge \|f\|_X\ge \|h\|_X.$ The remainder
of the argument mirrors that of Lemma 5.7.

Let us also sketch the argument when $p=\infty$ (i.e. for
Theorem 5.4).  Analogously to Lemma 5.7 we note that
$\|g\|_X \ge \|f\|_X \ge \|h\|_X$ where $g=\sum_{k\in
A_0}\|f_{I_k}\|_{\infty}e_{d_{k+1}}$ and $h=\sum_{k\in
A_0}\|f_{I_k}\|_{\infty}e_{d_k}.$ The remainder of the
argument is the same.\bull\enddemo

Now let us turn to the proofs of Theorems 5.4 and 5.5.
Suppose first that the couple $(X,Y)$ is a Calder\'on
couple.  Then the couple $(E_X(J), E_Y(J))$ must also be a
Calder\'on couple since there is a common averaging
projection from $(X,Y)$ onto $(E_X,E_Y)$.  Now it is clear
that $(E_X(J),E_Y(J))$ is exponentially separated (when $J$
is indexed as a sequence).  We can thus apply Theorem 4.2 to
obtain that $E_Y(J)$ has (LSP) and $E_X(J)$ has (LSP).  We
conclude this direction of the proof by showing that if
$E_Y(J)$ (and hence $F_Y$) has (LSP) then $E_Y$ has (LSP)
and so $Y$ is compressible.  A very similar argument shows
that $X$ is stretchable.

To prove this we suppose that $\{f_j,g_j\}_{j\in B}$ is an
interlaced pair of positive sequences in $E_Y$ with
$\|f_j\|_Y\le \|g_j\|_Y=1.$ For given $j$ let $l(j)$ be the
largest $k$ so that $P_{I_k}f_j\neq 0.$ (Note here if such a
largest $k$ does not exist then $j$ is the maximal element
of $B$ and $g_j=0;$ hence this case can be ignored.)  We
then split $f_j=f_j'+f_j''$ where $f_j''=P_{I_{l(j)}}f_j.$
Similarly we let $g_j=g_j'+g_j''$ where
$g_j'=P_{I_{l(j)}}g_j.$ Let $B_0=\{j:\|f_j'\|_Y\ge 1/2\},$
and let $B_1=B\setminus B_0.$

For $j\in B_0$ we set $v_j=(\|P_{I_k}f_j'\|_p)_{k\in A}\in
F_Y$; for $j\in B_1$ we set $v_j=(\|P_{I_k}f_j''\|_p)_{k\in
A}.$ For all $j\in B$ we set $w_j'=(\|P_{I_k}g_j'\|_p)$ and
$w_j''=(\|P_{I_k}g_j''\|_p).$

Let $(\alpha_j)_{j\in B}$ be positive and finitely nonzero.
First observe that $j$ in $B_0$ we must have supp $v_j<$
supp $(w_j'+w_j'')$.  Further $\|v_j\|_{F_Y}\ge (2C_0)^{-1}$
while $\|w_j+w_j''\|_{F_Y}\le C_0.$ Thus, since $F_Y$ has
(LSP) applying Lemma 2.6 we get the existence of a constant
$C_1$ depending only on $(E,F)$ so that $$ \|\sum_{j\in
B_0}\alpha_j(w_j'+w_j'')\|_{F_Y}\le C_1\|\sum_{j\in
B_0}\alpha_jv_j\|_{F_Y}.$$ Notice also that
$(w_j'+w_j'')_{j\in B_0}$ have disjoint supports so that we
can conclude that $$ \|\sum_{j\in B_0}\alpha_jg_j\|_Y \le
C_0 \|\sum_{j\in B_0}\alpha_j(w_j'+w_j'')\|_{F_Y}.$$
Similarly $$ \|\sum_{j\in B_0}\alpha_jv_j\|_{F_Y}\le
C_0\|\sum_{j\in B_0}\alpha_jf_j\|_Y.$$ Combining we have $$
\|\sum_{j\in B_0}\alpha_jg_j\|_Y\le C_0^2C_1 \|\sum_{j\in
B_0}\alpha_jf_j\|_Y.$$

We now obtain a similar estimate on $B_1.$ In fact, if we
set $B_2=\{j\in B_1:w_j''\neq 0\}$ then we can argue as
above to show that $$ \|\sum_{j\in
B_2}\alpha_jw_j''\|_{F_Y}\le C_1\|\sum_{j\in
B_2}\alpha_jv_j\|_{F_Y}$$ and hence obtain an estimate $$
\|\sum_{j\in B_1}\alpha_jg_j''\|_Y \le C_0^2C_1\|\sum_{j\in
B_1}\alpha_jf_j\|_Y.$$

Finally we observe that for $j\in B_1,$ $\|P_{l(j)}g_j\|_p
\le 4\|P_{l(j)}f_j\|_p$ by Lemma 5.6.  Thus for any $k$ $$
\align \|P_{I_k}\sum_{j\in B_1}\alpha_jg_j'\|_p &=(
\sum_{l(j)=k}|\alpha_j|^p\|P_{I_k}g_j'\|_p^p)^{1/p}\\ &\le
4(\sum_{l(j)=k}|\alpha_j|^p\|P_{I_k}f_j''\|_p^p)^{1/p}\\ &=
4 \|P_{I_k}\sum_{j\in B_1}\alpha_jf_j''\|_p.  \endalign $$
Thus $$ \|\sum_{j\in B_1}\alpha_jg_j'\|_Y \le
4C_0^2\|\sum_{j\in B_1}\alpha_jf_j''\|_Y.$$

Combining these estimates gives that $$ \|\sum_{j\in
B}\alpha_jg_j\|_{Y}\le C\|\sum_{j\in B}\alpha_jf_j\|_{Y}$$
for a suitable constant $C.$ This completes the proof that
$Y$ is compressible and, as explained above a similar
argument shows that $X$ is stretchable.

We now consider the other direction in Theorems 5.4 and 5.5.
We suppose $X$ is stretchable and $Y$ is compressible.  It
follows that $E_X$ has (RSP) and $E_Y$ has (LSP) and we can
apply both Lemmas 5.7 and 5.8.  We can immediately deduce:

\proclaim{Lemma 5.9}There exists $C$ so that if $0\le f,g\in
E_X+E_Y$ and $\|f_{I_k}\|_p \ge \|g_{I_k}\|_p$ for all $k\in
A$ then there exists $0\le T\in \Cal A(E_X,E_Y)$ with
$\|T\|_{(E_X,E_Y)}\le C$ and $Tf=g.$ \endproclaim

Now suppose $f,g\ge 0$ in $E_X+E_Y$ and that $K(t,g)\le
K(t,f)$ for all $t\ge 0.$ We define $f'=\sum_{k\in
A}2^{-d_k/p}\|f_{I_k}\|_pe_{d_k}$ and $g'=\sum_{k\in
A}2^{-d_k/p}\|g_{I_k}\|_pe_{d_k}.$ Then Lemma 5.9 yields the
conclusion that $K(t,g')\le CK(t,g)\le K(t,f)\le
C^2K(t,f').$ Now $(E_X(J),E_Y(J))$ is exponentially
separated.

Now for Theorem 5.5 we quote Theorem 4.5 to give that
$(E_X(J),E_Y(J))$ is a Calder\'on couple and hence there
exists $S\in \Cal A(E_X(J),E_Y(J))$ with
$\|S\|_{(E_X(J),E_Y(J))}\le C_2$, where $C_2$ depends only
on $(E,F)$, and $Sf'=g'.$ It follows easily from Lemma 5.9
that $(E_X,E_Y)$ and hence $(X,Y)$ is a uniform Calder\'on
couple.

In the case of Theorem 5.4 we note that it suffices to
consider the case when $f$ and $g$ are decreasing functions;
then $f'$ and $g'$ are also decreasing.  Then $K(t,g')\le
C^2K(t,f')$ for all $t$ implies that $$\|g'\chi_{[0,t]}\|_X
\le C^2\|f'\chi_{[0,t]}\|_X.$$ We further note that
$(E_X(J),\ell_{\infty}(J))$ has (RSP) by Lemma 3.4 and then
apply Lemma 4.3 to obtain a positive $S\in\Cal
A(E_X(J),\ell_{\infty}(J))$ with
$\|S\|_{E_X(J),\ell_{\infty}(J)}\le C_2$ and $Sf'=g'.$ This
leads to the desired conclusion.\bull\enddemo

\proclaim{Corollary 5.10}Let $X$ be an r.i. space on $[0,1]$
or $[0,\infty)$.  Suppose $X$ is $r$-concave for some
$r<\infty.$ In order that both $(L_1,X)$ and
$(L_{\infty},X)$ be Calder\'on couples it is necessary and
sufficient that $X$ be elastic.\endproclaim

\demo{Examples}We begin with the obvious remark that the
spaces $L_p$ for $1\le p\le\infty$ are elastic and so our
results include the classical results cited in the
introduction.  On the space $[0,\infty)$ one can basically
separate behavior at $\infty$ from behavior at $0$ so that
spaces of the form $L_p+L_q$ and $L_p\cap L_q$ are also
elastic.  Note however that we cannot apply Theorems 5.3 or
5.5 unless we have appropriate assumptions on either the
Boyd indices or convexity/concavity assumptions; thus pairs
of of such spaces are not always Calder\'on couples.

Let us now specialize to $[0,1].$ In certain special cases
we can easily see that an r.i. space is elastic.  For
example, suppose $X$ is the Lorentz space on $[0,1],$ for
which $q_X<\infty.$ Then it is immediately clear that $X$ is
elastic since $E_X$ is a weighted $\ell_p-$space.  Rather
more obscure elastic spaces can be built using a weighted
Tsirelson space for $E_X$.

On the other hand, it is possible to give easy examples
where $E_X$ fails (RSP) or (LSP).  Indeed if one takes any
symmetric sequence space $E$ on $\bold J$ which is not an
$\ell_p$-space and considers $E(w^n)$ where $1<w<2$ then
there is an r.i. space $X$ for which $E_X=E(w^n)$.  By
Proposition 2.3 $E_X$ fails (RSP) and (LSP).  In this case
we note that since $\kappa_+(E_X)=w$ and
$\kappa_-(E_X)=w^{-1},$ we have $p_X=q_X=(\log_2w)^{-1}.$ If
say $E=\ell_F(\bold Z_-)$ for some Orlicz function $F$
satisfying the $\Delta_2-$condition then $X$ is an
``Orlicz-Lorentz space'' given by $$ \|f\|_X \sim \int_0^1
F(f^*(t)t^{-1/p})\frac{dt}{t}$$ where $p=p_X=q_X.$ Note that
for such a space the pair $(L_{\infty},X)$ fails to be a
Calder\'on couple.  This answers a well-known question (cf.
[8],[28]).

In the next section we will investigate Orlicz spaces in
more detail.  We will also give examples of Orlicz spaces
$L_F$ for which $(L_{\infty},L_F)$ is not a Calder\'on
couple.

We will conclude this section by considering a situation
suggested by the example of Ovchinnikov [34] (cf.  [29]).

\proclaim{Theorem 5.11}Suppose $1<p<\infty$ and that $X$ is
an r.i. space on $[0,\infty)$ whose Boyd indices satisfy
either $q_X<p$ or $p<p_X\le q_X<\infty$.  Then $(X\cap
L_p,X+L_p)$ is a Calder\'on couple if and only if $X$ is a
Lorentz space of order $p.$\endproclaim

\demo{Proof}If $X$ is a Lorentz space of order $p,$ then
both $X+L_p$ and $X\cap L_p$ are also Lorentz spaces of
order $p,$ and so form a Calder\'on couple.  Conversely
suppose $(X\cap L_p,X+L_p)$ is a Calder\'on couple; then so
is $(E_{X\cap L_p},E_{X+L_p})$.  Let us consider the case
$q_X<p;$ the other case is similar.  Then $E_0=E_{X\cap
L_p}=E_X(\bold Z_-)\oplus E_{L_p}(\bold Z_+)$ and
$E_1=E_{X+L_p}=E_{L_p}(\bold Z_-)\oplus E_X(\bold Z_+).$
Note that for all $n$ we have $\|e_n\|_{X\cap L_p}\ge
\|e_n\|_{X+L_p}$; further if we rearrange the sequence
$(e_n)_{n\in\bold Z}$ so that $\|e_n\|_{X\cap
L_p}/\|e_n\|_{X+ L_p}$ increases, it is not difficult to see
that $(E_0,E_1)$ is exponentially separated.  Thus $E_0$ has
have (RSP) and $E_1$ has (LSP) for this ordering.  It also
follows easily form our assumptions on the Boyd indices that
there exists $k$ so that the gap in the new ordering for
$E_0$ between two consecutive elements of $\bold Z_+$ is at
most $k$.  Indeed the ratio $\|e_n\|_{X\cap
L_p}/\|e_n\|_{X+L_p}$ behaves like $2^{-n/p}\|e_n\|_X$ for
$n<0$ and like $2^{n/p}\|e_n\|_X^{-1}$ for $n\ge 0$ and we
have an estimate for $k>0,$ $C^{-1}2^{k/r} \le
\|e_{n+k}\|_X/\|e_n\|_X\le C2^k$ for suitable $C$ and $r$
with $q_X<r<p.$ Thus $E_0$ must be a weighted $\ell_p$-space
by the argument of Proposition 2.3.  It follows that
$E_X(\bold Z_-)$ is a weighted $\ell_p-$space.  Similarly
$E_X(\bold Z_+)$ is a weighted $\ell_p-$space and so $X$ is
a Lorentz space of order $p.$\bull\enddemo

\vskip2truecm

\subheading{6.  Orlicz spaces}

Let $F$ be an Orlicz function, i.e. a strictly increasing
convex function $F:[0,\infty)\to [0,\infty)$ satisfying
$F(0)=0.$ We will also assume that $F$ satisfies the
$\Delta_2-$condition with constant $\Delta$ i.e.  $F(2x)\le
\Delta F(x)$ for every $x>0.$ We will use the notation
$F_t(x)=F(tx)/F(x).$

We recall first that $F$ is said to be {\it regularly
varying at $\infty$ (resp. at $0$)}, in the sense of
Karamata, if the limit $\lim_{t\to\infty} F_t(x)$ (resp.
$\lim_{t\to 0}F_t(x)$) exists for all $x$ (in fact, it
suffices that the limit exists when $x\le 1.$) In this case
there exists $p, \ 1\le p<\infty$ so that
$\lim_{t\to\infty}F_t(x)=x^p$ (resp.  $\lim_{t\to
0}F_t(x)=x^p$); $F$ is then said to be regularly varying
with order $p$.  See [6] for details.

\proclaim{Lemma 6.1}The following conditions are
equivalent:\newline (1) $F$ is equivalent to an Orlicz
function $G$ which is regularly varying with order $p$ at
$\infty$ (resp.  $0$).\newline (2) There exists a constant
$C$ so that if $x_0\le 1$ there exists $0<t_0<\infty$ so
that if $t\ge t_0$ (resp.  $t\le t_0$) and $x_0\le x\le 1,$
$$ C^{-1}x^p \le F_t(x) \le Cx^p.$$ (3) There exists a
constant $C$ so that if $x\le 1$ $\limsup_{t\to\infty}F_t(x)
\le C\liminf_{t\to\infty}F_t(x).$ (resp.  $\limsup_{t\to
0}F_t(x)\le C\liminf_{t\to 0}F_t(x).$) \endproclaim

\demo{Proof}The implication $(1)\Rightarrow (3)$ is
immediate and $(3)\Rightarrow (2)$ is a simple compactness
argument.  We indicate the details of $(2)\Rightarrow (1).$
Let $f(x)=\log F(e^x)$ for $x\in\bold R.$ The function
$f(x)-x$ is then increasing.  Then it is easy to translate
(2) as:

\flushpar $(2)':$ there exists $c$ so that if $y_0\ge 0$
there exists $x_0$ so that if $0\le y\le y_0,$ then $$
|f(x)-f(x-y)-py| \le c,$$ whenever $x\ge x_0.$

Now we can pick a function $u=u(x)$ for $x\in\bold R$ so
that $u(x)=0$ for $x\le 0$, $u$ is differentiable,
increasing, $u'(x)\le 1,$ $\lim_{x\to\infty}u(x)=\infty,$
and $|f(x)-f(x-y)-py|\le c$ for $0\le y\le u(x).$ Now define
$g(x)=f(x)$ if $u(x)=0$ and $$ g(x)= \frac1u\int_{x-u}^x
(f(s)+p(x-s))ds$$ if $u>0.$ It is easy to show that $f-g$ is
bounded.  Further if $u>0,$ $$ g'(x) =
-\frac{u'}{u}(g(x)-f(x-u)-pu) + \frac1u(f(x)-f(x-u)-pu) +p.
$$ Since $$ \frac{u'(x)}{u(x)} \le \frac1u$$ it is easy to
see that $\lim g'(x)=p$ and so if $G_0(x)=\exp g(\log x)$
then $G_0$ is regularly varying and equivalent to $F$.

It remains to construct a convex $G$ with the same
properties.  First note that since $f(x)-x$ is increasing we
have if $u>0,$ $$ g(x)-f(x) \le \frac1u \int_{x-u}^x
(p-1)(x-s)ds\le \frac{p-1}2u.$$ Hence $$ \align g'(x) &\ge
p+\frac{1-u'}{u}(f(x)-f(x-u)-pu) -\frac{p-1}2u'\\ &\ge p-
(p-1)(1-u') -(p-1)u'\\ &\ge 1. \endalign $$ It now follows
that $G_0(x)/x$ is increasing.  The proof is completed by
setting $G(x)=\int_0^x G_0(x)/x dx$ and it is then easy to
verify that $G$ has the desired properties.  \bull\enddemo

If $F$ is an Orlicz function, $0<x\le 1,$ and $C>1$ we can
define $\Psi_p^{\infty}(x,C)$ (resp.  $\Psi_p^{0}(x,C)$) to
be the supremum (possibly $\infty$) of all $N$ so that there
exist $a_1<a_2<\cdots<a_N$ with $a_k/a_{k-1}\ge 2$ for $k\le
N-1$ and $a_1\ge 1$ (resp.  $a_N\le 1$) so that for all $k$
either $F_{a_k}(x) \ge Cx^p$ or $x^p \ge CF_{a_k}(x).$ It is
easy to show that

\proclaim{Proposition 6.2}$F$ is equivalent to a regularly
varying function of order $p$ at $\infty$ (resp. at $0$) if
and only if for some $C$ and all $0<x\le 1$ we have
$\Psi_p^{\infty}(x,C)<\infty$ (resp.
$\Psi_p^{0}(x,C)<\infty).$\endproclaim

We omit the proof which is immediate.  However we can now
state the result of Montgomery-Smith [33] which
characterizes Orlicz spaces which are Lorentz spaces (see
Lorentz [26]).

\proclaim{Theorem 6.3}In order that $L_F[0,1]$ coincides
with a Lorentz space of order $p$ it is necessary and
sufficient that there exist $C_0,C_1$ and $r>0$ so that for
every $x$ with $0<x\le 1$ we have $\Psi_p(x,C_0)\le
C_1x^{-r}.$\endproclaim

This is a somewhat disguised restatement of
Montgomery-Smith's result.  However we will not pause in our
exposition to derive this result as a proof is implicit in
our approach to elastic Orlicz spaces.  Further, we state
the result in order to motivate the following definition.

For $C>1$ and $0<x\le 1$ let us define
$\Phi_+^{\infty}(x,C)$ (resp.  $\Phi_+^{0}(x,C)$ to be the
supremum of all $n$ so that there exist $a_1<b_1\le
a_2<b_2\cdots\le a_n<b_n$ with $a_1\ge 1$ (resp.  $b_n\le
1$) so that $F_{b_k}(x) \ge CF_{a_k}(x),$ for $1\le k\le n.$
For $C>1$ and $0<x\le 1$ let us define
$\Phi_-^{\infty}(x,C)$ (resp.  $\Phi_-^{0}(x,C)$) to be the
supremum of all $n$ so that there exist $a_1<b_1\le
a_2<b_2\cdots\le a_n<b_n$ with $a_1\ge 1$ (resp.  $b_n\le
1$) so that $F_{a_k}(x) \ge CF_{b_k}(x),$ for $1\le k\le n.$
We say that $F$ is {\it elastic at $\infty$ (resp. at $0$)}
if there exist $C_0,C_1>1$ and $r>0$ so that for $0< x\le 1$
we have $\Phi_+^{\infty}(x,C_0)+\Phi_-^{\infty}(x,C_0)\le
C_1x^{-r}$ (resp.  $\Phi_+^0(x,C_0)+\Phi_-^0(x,C_0) \le
C_1x^{-r}$).  From now on, we will consider only the case at
$\infty$ although similar results can always be proved at
$0.$

\proclaim{Lemma 6.3}$F$ is elastic at $\infty$ if and only
if there exist constants $C_0,C_1>1$ and $r>0$ so that if
$0<x\le 1$, $\Phi_+^{\infty}(x,C_0)\le C_1x^{-r}$ (resp.
$\Phi_-^{\infty}(x,C_0)\le C_1x^{-r}$).\endproclaim

\demo{Proof}Assume $\Phi_+^{\infty}(x,C_0)\le C_1x^{-r}.$
Suppose $1\le a_1<b_1\le\cdots\le a_n<b_n$ with $F_{a_k}(x)
\ge eC_0 F_{b_k}(x)$ for $1\le k\le n$.  Consider an
interval $[b_k,a_{k+1}]$ where $1\le k\le n-1.$ Let
$\nu=\nu_k$ be the integer part of $(\log C_0)^{-1}(\log
F_{a_{k+1}}(x) -\log F_{b_k}(x)).$ Then we can find
$b_k=c_0<c_1<\cdots<c_\nu\le a_{k+1}$ so that $\log
F_{c_k}(x)-\log F_{c_{k-1}}(x)=\log C_0.$ It follows that $$
\sum_{k=1}^{n-1}\nu_k \le \Phi_+(x,C_0)$$ and hence that $$
\sum_{k=1}^{n-1}(\log F_{a_{k+1}}(x)-\log F_{b_k}(x)) \le
(\log C_0)(\Phi_+(x,C_0)+n-1) $$ and thus $$ \log
F_{b_n}(x)-\log F_{a_1}(x) \le (\log C_0)(C_1x^{-r}-1)-n.$$
Now $$\log F_{b_n}(x)-\log F_{a_1}(x) \ge \log F_{b_n}(x)
\ge -C_2|\log x|-C_3$$ for suitable $C_2,C_3$ by the
$\Delta_2$ condition.  Hence $$ n \le (\log C_0)(C_1
x^{-r}-1) + C_2|\log x| +C_3$$ and so $$
\Phi_-^{\infty}(x,eC_0) \le C_4x^{-2r}$$ for a suitable
$C_4.$ The other case is similar.\bull\enddemo

\proclaim{Proposition 6.4}The following conditions on $F$
are equivalent:  \newline (1) $F$ is elastic at
$\infty$.\newline (2) There exist constants $C_0,C_1>1$ so
that if $1\le a_1<b_1\le\cdots \le a_n<b_n$ and $0\le x\le
1$ then:  $$ \sum_{k=1}^n (F_{b_k}(x)-C_0F_{a_k}(x)) \le
C_1.$$ (3) There exist constants $C_0,C_1>1$ so that if
$1\le a_1<b_1\le\cdots \le a_n<b_n$ and $0\le x\le 1$ then:
$$ \sum_{k=1}^n (F_{a_k}(x)-C_0F_{b_k}(x)) \le C_1.$$ (4)
There exists a bounded monotone increasing function
$w:[1,\infty)\to\bold R$ and a constant $C_0$ so that if
$1\le s\le t$ and $0\le x\le 1$ then $$ F_t(x) \le CF_s(x)
+w(t)-w(s).$$ (5) There exists a bounded monotone increasing
function $w:[1,\infty)\to\bold R$ and a constant $C_0$ so
that if $1\le s\le t$ and $0\le x\le 1$ then $$ F_s(x) \le
CF_t(x) +w(t)-w(s).$$ \endproclaim

\demo{Proof}$(1)\Rightarrow (2).$ We assume that for
suitable constants $C_2,C_3>1$ and $r>0,$ we have
$$\Phi_+^{\infty}(x,C_2)\le C_3x^{-r}.$$ We will assume that
$C_2>\Delta$ from which it follows easily
$F_{b}(x)/F_{a}(x)\ge C_2$ implies that $b>2a.$ First
suppose $m$ is an integer with $m>r$.  We will estimate
$\Phi_+^{\infty}(x,C_2^m).$ Suppose $1\le
a_1<b_1\le\cdots\le a_n<b_n$ and
$F_{b_k}(x)>C_2^mF_{a_k}(x).$ Let $s$ be the smallest
integer greater than $|\log_2x|+1.$ Then $a_s \ge x^{-1}$
and $a_{ks}>x^{-1}b_{(k-1)s}$ for $2\le k\le [n/s].$ Let
$\xi=x^{1/m}.$

Now for each $1\le k\le [n/s]$ there exists $\sigma_k$ with
$0\le \sigma_k\le m-1$ so that
$F_{\xi^{\sigma_k}b_{sk}}(\xi)\ge C_2
F_{\xi^{\sigma_k}a_{sk}}(\xi)$ and the intervals
$[\xi^{\sigma_k}a_{sk},\xi^{\sigma_k}b_{sk}]$ are disjoint
in $[1,\infty).$ Hence we have an estimate that
$$\Phi_+^{\infty}(\xi,C_2) \ge [n/s]$$ and this means that
$$ [n/s] \le C_3\xi^{-r}.$$ Thus $$ n \le C_3(s+1)x^{-r/m}
\le C_4+C_5|\log x|x^{-r/m}$$ for suitable constants
$C_4,C_5.$ This leads to an estimate
$$\Phi_+^{\infty}(x,C_2^m)\le C_6x^{-\alpha}$$ where $0<
\alpha<1.$

Now suppose $C_0=\Delta C_2^m.$ Suppose $1\le
a_1<b_1\le\cdots\le a_n<b_n$ and that $0\le x_k\le 1$ for
$1\le k\le n.$ For $j\in\bold N$ let $I_j$ be the set of $k$
such that $2^{-j}< x_k\le 2.2^{-j}.$ Then $$ \align
\sum_{k\in I_j}(F_{b_k}(x_k)-C_0F_{a_k}(x_k)) &\le
\sum_{k\in I_j}(\Delta F_{b_k}(2^{-j})-C_0F_{a_k}(2^{-j}))\\
&\le \Delta \Phi_+^{\infty}(2^{-j},C_2^m)\max_k
F_{b_k}(2^{-j})\\ &\le C_6\Delta 2^{-(1-\alpha) j}.
\endalign $$ Thus $$ \sum_{k=1}^n(F_{b_k}(x)-C_0F_{a_k}(x))
\le C_6\Delta\sum_{j=1}^{\infty} 2^{-(1-\alpha) j}.$$ This
establishes (2).

$(2)\Rightarrow (4).$ We define $w(t)$ for $1\le t<\infty$
by setting $w(t)$ to be the supremum of
$\sum_{k=1}^n(F_{b_k}(x_k)-C_0F_{a_k}(x_k))$ over all $n$
and all $1\le a_1<b_1\le\cdots\le a_n<b_n\le t$ and all
$0\le x_k\le 1$ for $1\le k\le n.$ Clearly $w(t)$ is
increasing and bounded above by $C_1$.  Condition (4) is
immediate from the definition.

$(4)\Rightarrow (1).$ Suppose $0<x\le 1$ and that $1\le
a_1<b_1\le\cdots\le a_n<b_n,$ are such that
$F_{b_k}(x)>2C_0F_{a_k}(x).$ Then we have $$ C_0\sum_{k=1}^n
F_{b_k}(x) \le \sum_{k=1}^n(w(b_k)-w(a_k))\le C_1$$ where
$C_1=\lim_{x\to\infty}w(x) -w(1).$ Now $F_{t}(x)\ge C_2x^r$
for all $t$, for a suitable $C_2$, by the
$\Delta_2-$condition.  Thus $$ \Phi_+^{\infty}(x,2C_0) \le
C_1(C_0C_2)^{-1}x^{-r}.$$ The implication now follows from
Lemma 6.3.

The remaining implications are similar.\bull\enddemo

\proclaim{Lemma 6.5}If $F$ is elastic at $\infty$ then $F$
is equivalent to an Orlicz function which is regularly
varying at $\infty.$\endproclaim

\demo{Proof}It follows immediately from (4) above that
$$\limsup_{t\to\infty}F_t(x) \le
C_0\liminf_{t\to\infty}F_t(x)$$ for $0<x\le 1$.  Apply Lemma
6.1.\bull\enddemo

We now come to our main theorem on elastic Orlicz functions.

\proclaim{Theorem 6.6}Let $F$ be an Orlicz function
satisfying the $\Delta_2-$condition.  Then the following are
equivalent:\newline (1) $F$ is elastic at $\infty.$\newline
(2) $L_F[0,1]$ is stretchable.\newline (3) $L_F[0,1]$ is
compressible.\newline (4) $L_F[0,1]$ is elastic.
\endproclaim

\demo{Proof}We will only show $(1)\Rightarrow (2)$ and
$(2)\Rightarrow (1).$ The other implications will then be
clear.  We will write $E=E_F$ for $E_X$ where $X=L_F[0,1].$
Then $E_F$ is the modular sequence space of $\bold Z_-$
defined by $ \|x\|_{E_F} =1$ if and only if $\sum_{n\in\bold
Z_-}F(x(n))2^n=1.$ Let us define $\lambda_n$ for $n\in\bold
Z_-$ by $F(\lambda_n)=2^{-n}.$ Then $(\lambda_n)_{n\in \bold
Z_-}$ is strictly decreasing and $\lambda_{n-1}\le
2\lambda_{n}$ for $n<0.$

$(1)\Rightarrow (2).$ We must show that $E$ has (RSP).  To
show this it suffices to show the existence of a constant
$C$ so that if $a_1< b_1< c_1\le a_2<b_2<c_2\le\cdots\le
a_n<b_n<c_n\le 0$ and supp $x_k\subset [a_k,b_k)$, supp
$y_k\subset [b_k,c_k)$ and $\|y_k\|_E\le \|x_k\|_E=1$ then
$$ \|\sum_{k=1}^n\alpha_ky_k\|_E \le
C\|\sum_{k=1}^n\alpha_kx_k\|_E.$$

To do this let us suppose $n,a_k,b_k,c_k,x_k$ as fixed and
let $\Gamma$ be the least constant $C$ for which this
inequality holds.  We show a uniform bound on $\Gamma.$ We
can suppose the existence of constants $C_0,C_1$ and an
increasing function $w:[1,\infty)\to \bold R$ with
$\lim_{x\to\infty}w(x)=w(1)+C_1$ so that if $1\le s\le t$ $$
F_s(x)\le C_0F_t(x) +w(t)-w(s)$$ for $0\le x\le 1.$

Let us define $x_k'=\sum_{|x_k(j)|\ge
\frac12\lambda_{b_k}}x_k(j)e_j$ and also
$y_k'=\sum_{|y_k(j)|\ge \frac12\lambda_{c_k}}y_k(j)e_j.$

Then $$ \sum 2^jF(2|x_k(j)-x'_k(j)|) \le
2^{-b_k}\sum_{j<b_k}2^{j}\le 1.$$ Thus $\|x_k-x_k'\|_E\le
1/2$ and similarly $\|y_k-y_k'\|_E\le 1/2.$ We let
$x_k''(j)=\min(2|x_k'(j)|,\lambda_{a_k})$ and $y_k''(j)
=\min(2|y_k'(j)|,\lambda_{b_k}).$ Then
$\|x_k''\|_E,\|y_k''\|_E\le 2.$

Now for any $\alpha_1,\ldots,\alpha_n$ such that
$\|\sum_{k=1}^n\alpha_kx_k\|_E=1$ we set
$z=\sum_{k=1}^n\alpha_ky_k'$ and
$v=\sum_{k=1}^n\alpha_kx_k'.$ We also let $u=\sum_{k=1}^n
\alpha_k\lambda_{b_k}e_{b_k}.$ Then for fixed $k$, $$ \align
\sum_{j\in [b_k,c_k)}2^jF(|\alpha_ky_k'(j)|) &\le \sum 2^j
F(|\alpha_ky_k''(j)|)\\ &\le C_0
2^{b_k}F(|\alpha_k|\lambda_{b_k}) \sum_{y_k''(j)\neq 0 }
2^jF(y_k''(j))\\ &\hskip.5truecm +\sum_{y_k''(j)\neq
0}2^jF(y_k''(j))(w(\lambda_{b_k})-w(y_k''(j))) \\ &\le
C_0\Delta 2^{b_k}F(|\alpha_k|\lambda_{b_k}) +
\Delta(w(\lambda_{b_k})-w(\lambda_{c_k})).  \endalign $$

On summing, we get $$\sum_j 2^jF(|z(j)|) \le
C_0\Delta\sum_j2^jF(|u(j)|) + \Delta C_1.$$

Now in the other direction, for fixed $k$, $$
2^{b_k}F(|\alpha_k|\lambda_{b_k}) \le
C_0F(|\alpha_k|x_k''(j))F(x_k''(j))^{-1} +
w(\lambda_{a_k})-w(\lambda_{b_k})$$ whenever $x_k''(j)\neq
0.$

Thus $$ \align
2^{b_k}F(|\alpha_k|\lambda_{b_k})(\sum_j2^jF(x_k''(j)))\le &
C_0\sum_j2^jF(|\alpha_kx_k''(j)|) + \\
&+(\sum_j2^jF(x_k''(j)))(w(\lambda_{a_k}-w(\lambda_{b_k})).
\endalign $$

Now we observe that $1/2\le \|x_k'\|_E \le 1$ so that
$1/2\le \|x_k''\|_E\le 2$.  Hence $1/2 \le \sum_j
2^jF(x_k''(j)) \le \Delta.$ Thus we have:  $$
2^{b_k}F(|\alpha_k|\lambda_{b_k}) \le
2C_0\Delta\sum_j2^jF(|\alpha_kx_k'(j)|) + 2\Delta
(w(\lambda_{a_k})-w(\lambda_{b_k})).$$ Summing as before $$
\sum_j2^jF(|u(j)|) \le 2C_0\Delta\sum_j 2^jF(|v(j)|)
+2C_1\Delta.$$ We thus have an estimate $$ \sum_j2^j
F(|z(j)|)\le C_2\sum_j2^jF(|v(j)|) + C_3$$ for constants
$C_2,C_3$ depending only on $F.$ This in turn implies an
estimate $\|z\|_E \le C_4\|v\|_E$ for some constant $C_4.$

Now we conclude by noting that:  $$ \align
\|\sum_{k=1}^n\alpha_ky_k\|_E &\le \|z\|_E
+\|\sum_{j=1}^k\alpha_k(y_k-y_k')\|_E\\ &\le C_4\|v\|_E +
\Gamma\max_j\|y_j-y_j'\|_E \|\sum_{k=1}^n\alpha_kx_k\|_E\\
&\le (C_4+\frac{\Gamma}{2}) \endalign $$ Thus $\Gamma\le
C_4+\Gamma/2$ and so $\Gamma\le 2C_4$ and $E$ has (RSP).

$(2)\Rightarrow (1).$ Suppose $E$ has (RSP).  This implies
that for some $C_0,$ if $a_1<b_1<a_2<\cdots<a_n<b_n$ are
negative integers, and $0\le x_k $ with $ \sum
2^{a_k}F(\lambda_{a_k}x_k)\le 1$ then $\sum 2^{b_k}
F(\lambda_{b_k}x_k)\le C_0.$ We also note from the
$\Delta_2$-condition that we can suppose $C_1F_t(x)\ge x^r$
for some $C_1$ and $r$ and all $t>0,\ 0\le x\le 1,$

For any constant $C>2C_0\Delta^3$ and $0<x\le 1$ suppose now
that $1\le c_1<d_1\le c_2<\cdots<d_{n-1}\le c_n<d_n$ and
$F_{c_k}(x)\ge CF_{d_k}(x).$ Then we must have $d_k>2^3c_k.$
Now choose $b_{n-k+1}\in \bold Z_-$ to be the largest
integer so that $\lambda_{b_{n-k+1}}>c_k$ and let
$a_{n-k+1}$ be the smallest integer so that
$\lambda_{a_{n-k+1}}<d_k.$ It is clear that $a_1
<b_1<a_2\cdots<a_n<b_n.$ Further $\lambda_{b_{n-k+1}}\le
2c_k$ and $\lambda_{a_{n-k+1}} \ge d_k/2.$ It follows that
for every $k$ with $1\le k\le n$ we have
$$2^{b_k}F(\lambda_{b_k}x) \ge
C\Delta^{-2}2^{a_k}F(\lambda_{a_k}x).$$

Now suppose $n>C_1x^{-r}.$ Then we can select a subset $J$
of $\{1,2,\ldots,n\}$ so that $1/2 \le \sum_{k\in
J}2^{a_k}F(\lambda_{a_k}x)\le 1.$ Then we can conclude that
$$\frac{C\Delta^{-2}}2 \le\sum_{k\in
J}2^{b_k}F(\lambda_{b_k}x) \le C_0.$$ Since $C>2C_0\Delta^2$
we reach a contradiction and conclude that $n\le C_1x^{-r}.$
Thus $$\Phi_-^{\infty}(x,C)\le C_1x^{-r}$$ and $F$ is
elastic by Lemma 6.3.\bull\enddemo

Of course there are corresponding results for sequence
spaces and Orlicz spaces on $[0,\infty)$.  We will omit the
proofs.

\proclaim{Theorem 6.7}Suppose $F$ is an Orlicz function
satisfying the $\Delta_2-$condition.  Then:\newline (1) In
order that $\ell_F$ be elastic (resp. compressible, resp.
stretchable) it is necessary and sufficient that $F$ be
elastic at $0$.\newline (2) In order that $L_F[0,\infty)$ be
elastic (resp. compressible, resp. stretchable) it is
necessary and sufficient that $F$ be elastic at both $0$ and
$\infty.$\endproclaim

\demo{Remark}It is perhaps worth pointing out at this point
that the theorem of Montgomery-Smith (Theorem 6.3) cited
above can be proved in much the same manner as Theorem 6.6;
the problem in this case is to show that $E$ is a weighted
$\ell_p-$space.  In fact our proof of Theorem 6.6 is derived
from the arguments used by Montgomery-Smith [33].

Returning to the case of $[0,1]$ we note the following
simple deduction.

\proclaim{Proposition 6.8}If the Orlicz space $L_F[0,1]$ is
elastic then its Boyd indices $p_F=p_{L_F}$ and
$q_F=q_{L_F}$ coincide.\endproclaim

\demo{Proof}In fact we can suppose $F$ is regularly varying
by Lemma 6.5 and so the conclusion is
immediate.\bull\enddemo

\demo{Remark}The analogous result holds for sequence spaces,
but not for $L_F[0,\infty)$ where one must consider behavior
at both $0$ and $\infty$.  Thus $L_p\cap L_q$ is elastic for
any $p,q.$ Let us also mention at this point that
Proposition 6.8 allows us very easily to give examples of
Orlicz function spaces $L_F[0,1]$ so that $(L_{\infty},L_F)$
is not a Calder\'on couple by simply ensuring that $p_F\neq
q_F$.\enddemo

\demo{Examples}We now give two examples to separate the
concepts implicit in our discussion above.  We first
construct a regularly varying Orlicz function which is not
elastic.  To do this first suppose $(\xi_n)$ is a positive
sequence, bounded by one and tending monotonically to zero.
We define $\phi(x)=2$ if $x\le 1$ and then
$\phi(x)=2+(-1)^n\xi_n$ if $2^{n-1}<x\le 2^n.$ Define
$f(x)=\int_0^x \phi(t)dt$ and $F(x)=\exp(f(\log x)).$ Then
$F(x)/x$ is increasing and hence $F_1(x)=\int_0^xF(t)/t\,
dt$ is an Orlicz function equivalent to $F$.  Further $F$
and $F_1$ are regularly varying of order 2. It remains to
show that $F_1$ or equivalently $F$ is not elastic at
$\infty.$ Suppose $C>1$ and that $ 0<x\le 1$.  If
$2^{n-1}>\log x^{-1}$ then $$|\log
(\frac{F(e^{2^n}x)}{F(e^{2^n})})-\log
\frac{F(e^{2^{n+1}}x)}{F(e^{2^{n+1}})}| \ge
(\xi_n+\xi_{n+1})\log x^{-1}.$$ If we assume that $\xi_n$
goes to zero slowly enough, say $\xi_n\sim(\log\log n)^{-1}$
this will exceed $\log C,$ $O(\exp(x^{-r}))$ times for some
$r>0$ and so $F_1$ cannot be elastic.

Our second construction is of an elastic Orlicz space which
is not a Lorentz space.  It is of course clear that
conversely that every Lorentz space is elastic.  We note
first that if $F(x)=\exp(f(\log x))$ where $f$ is convex
then $F$ is elastic at $\infty$, by applying Proposition 6.4
(4) (equally the same conclusion holds when $f$ is concave).
We thus consider a function $\phi(t)=2+\psi(t)$ where
$\psi(t)$ is bounded by one and decreases monotonically to
$0.$ Let $f(x)=\int_0^t\phi(t)dt$ as above.  As usual it may
be necessary to convexify $F$ by constructing $F_1$; however
this is equivalent to $F.$ Now we show that for $L_F[0,1]$
to be a Lorentz space it is necessary that $\psi$ tends to
zero at a certain rate.  In fact if
$\Psi_2^{\infty}(x,C_0)\le C_1x^{-r}$ it follows that
$\psi(2^{C_1 x^{-r}+1}) < \log C_0/\log x^{-1}$ and hence
that $\psi(u)=O((\log\log u)^{-1}).$ Thus if we choose
$\psi$ converging to zero slowly enough then $L_F[0,1]$ is
an elastic non-Lorentz space.\enddemo

We now turn to the general problem of determining when a
pair of Orlicz spaces $L_F[0,1]$ and $L_G[0,1]$ forms a
Calder\'on couple.  Of course if the Boyd indices satisfy
$q_F<p_G$ this can only happen if both $F$ and $G$ are
elastic at $\infty$ in which case $p_F=q_F$ and $p_G=q_G.$
Brudnyi [8] has conjectured that if $L_F$ and $L_G$ are
distinct then if $(L_F,L_G)$ forms a Calder\'on couple then
we must have $p_F=q_F$ and $p_G=q_G.$ The next theorem shows
that that if either $p_F\neq q_F$ or $p_G\neq q_G$ then $F$
and $G$ must in some sense be similar functions.  However
following the theorem we will give a counterexample to
Brudnyi's conjecture.

\proclaim{Theorem 6.9}Suppose $F$ and $G$ are Orlicz
functions satisfying the $\Delta_2-$condition and such that
$(L_F[0,1],L_G[0,1])$ forms a Calder\'on pair.  Then either
$F$ and $G$ are both elastic or $p_F=p_G$ and
$q_F=q_G.$\endproclaim

\demo{Proof}Let us assume that $q_F>q_G.$ The other case is
similar.  It will be convenient to pick $q_0,q_1$ so that
$q_G<q_0<q_F<q_1$ and to suppose (by passing to equivalent
functions) that $F(x)/x^{q_1}$ and $G(x)/x^{q_0}$ are
decreasing.

Let $\Cal F$ be the closure of the set of functions
$\{F_t:t\ge 1\}$ in $C[0,1]$.  This set is relatively
compact.  For each $M>1$ let $\Cal F_M$ be the closure of
the set of functions $\{F_t:t\ge 1, \ F(t)/G(t)\ge M\}$ and
let $\Cal F_{\infty}=\cap_M \Cal F_M.$ Similarly if $a<1$ we
let $\Cal F_a$ be the closure of the set of functions
$\{F_t:t\ge 1,\ F(t)/G(t)\le a\}$ and set $\Cal F_0=\cap_a
\Cal F_a.$

Now suppose $t\ge 1$ and $A_t$ is a measurable subset of
$[0,1]$ such that $\mu(A_t)=F(t)^{-1}.$ Then
$\|\chi_{A_t}\|_{L_F}=t^{-1}$ while
$\|\chi_{A_t}\|_{L_G}=s^{-1}$ where $G(s)=F(t).$ If
$F(t)>G(t)$ we conclude that $s>t$ and further from the
$\Delta_2-$condition for $G$ we have that
$\|\chi_{A_t}\|_{L_G} \ge
\phi(F(t)/G(t))\|\chi_{A_t}\|_{L_F}$ where $\phi$ is a
function satisfying $\lim_{u\to\infty}\phi(u)=\infty.$

Suppose $\Cal F_{\infty}$ is nonempty and $H_1,H_2\in\Cal
F_{\infty}.$ Then we can find a sequence $(t_n)_{n\ge 1}$
such that $t_1\ge 2$, $t_n>2t_{n-1}$ and
$F(t_n)/G(t_n)\to\infty,$
$$\frac{\|\chi_{A_{t_n}}\|_{L_G}}{\|\chi_{A_{t_n}}\|_{L_F}}
\ge
2\frac{\|\chi_{A_{t_{n-1}}}\|_{L_G}}{\|\chi_{A_{t_{n-1}}}\|_{L_F}}$$
for $n\ge 1$ and such that $F_{t_{2n}}\to H_1$ while and
$F_{t_{2n-1}}\to H_2.$ Since $\mu A_{t_n}\le 2^{-n}$ we can
suppose these sets are disjoint.  If we restrict to the
sub$-\sigma-$algebra $\Cal A$ of the Borel sets generated by
$(A_{t_n})$ then $(L_F(\Cal A),L_G(\Cal A))$ forms a
Calder\'on couple.  Regarded as a couple of sequence spaces
it is exponentially separated and hence the Orlicz modular
space $\ell_{F_{t_n}}$ has (LSP) by Theorem 4.2.  By passing
to a subsequence of the unit vectors it follows that both
the Orlicz sequence spaces $\ell_{H_1}$ and $\ell_{H_2}$
have (LSP) and further that the space obtained by
interlacing their bases has (LSP).  Hence from Proposition
2.4 $H_1(x)$ and $H_2(x)$ are both equivalent to some
(common) $x^{p_0}.$ We thus conclude that there exists $p_0$
so that any $H\in\Cal F_{\infty}$ is equivalent to
$x^{p_0}.$

By similar reasoning, if $\Cal F_0$ is nonempty there exists
$p_1$ so that every $H\in\Cal F_0$ is equivalent to
$x^{p_1}.$

Now suppose $q_0<r_1<r_2<q_F.$ We pick $m$ an integer large
enough so that $(m-2)r_1+2q_1<mr_2.$ Then for any $\xi<1$
the function $F_0(x)=\max\{x^{r_2}t^{-r_2}F(t):\xi^mx\le
t\le x\}$ is equivalent to $F$ and therefore
$F_0(x)/x^{r_2}$ cannot be decreasing everywhere.  Thus for
any $x_0$ there exists $x\ge x_0$ such that for some
$\delta>0$ we have $F_0(u)/u^{r_2}<F_0(x)/x^{r_2}$ if
$x-\delta<u<x.$ It follows that $F_0(x)=F(x)$ and hence
$F(x) \ge x^{r_2}t^{-r_2}F(t)$ if $\xi^mx\le t\le x.$

Next define $F_1(y)=\max\{y^{r_1}t^{-r_1}F(t):\xi y\le t\le
y\}.$ Notice that $F_1(y) \le \xi^{q_1-r_1}F(\xi y).$ We
will argue that $F_1(x)/x^{r_1}$ cannot be decreasing on
$(\xi^{m-1}x,\xi x)$.  If it is then $F_1(\xi x)\le
F_1(\xi^{m-1}x)$ and hence $$ \align F(x) &\le
\xi^{-q_1}F(\xi x)\\ &\le \xi^{-q_1}F_1(\xi x)\\ &\le
\xi^{-(m-2)r_1-q_1}F_1(\xi^{m-1}x)\\ &\le
\xi^{-(m-2)r_1-2q_1}F(\xi^mx) \endalign $$ and hence
$(m-2)r_1 +2q_1>mr_2$ contrary to assumption.

We now argue as above and conclude similarly that there
exists $u$ with $\xi^{m-1}x\le u\le \xi x$ such that
$F(u)\ge u^{r_1}t^{-r_1}F(t)$ for $\xi u\le t\le u.$

Now notice that $$F(u)/G(u) \le
u^{r_2-q_0}x^{q_0-r_2}F(x)/G(x)$$ and so $$F(u)/G(u) \le
\xi^{r_2-q_0}F(x)/G(x).$$

It follows that given any $\xi<1$ any $x_0$ we can pick
$x\ge x_0$ so that $F(x)\ge x^{r_1}t^{-r_1}F(t)$ for $\xi
x\le t\le x$ and either $F(x)/G(x) \ge \xi^{-(r_2-q_0)/2}$
or $F(x)/G(x) \le \xi^{(q_0-r_2)/2}.$ Thus we can find a
sequence $t_n\to\infty$ such that $F_{t_n}(x)\ge x^{r_1}$
for $n^{-1}\le x\le 1$ and either $F(t_n)/G(t_n)\to\infty$
or $F(t_n)/G(t_n)\to 0.$

Consider the former case.  Then there exists $H\in \Cal
F_{\infty}$ with $H(x)\ge x^{r_1}.$ Hence $p_0 \ge r_1.$ In
the latter case $p_1\ge r_1.$ Since $r_1<q_F$ is arbitrary
we conclude that either $p_0=q_F$ or $p_1=q_F.$

Consider the case $p_0=q_F$; in particular $F(t)/G(t)$ is
unbounded for $t\ge 1.$ We will argue that $F(t)/G(t)$ tends
to infinity.  For any $\xi<1$ consider the function
$h(x)=\min\{F(t)/G(t):\xi x\le t\le x\}.$ If $h$ does not
converge to $\infty$ then given any $M$ and $x_0$ there
exists $x>x_0$ and $\delta>0$ so that $h(x)=M$ and
$h(x)<h(u)$ for $x-\delta<u<x$ and this implies that
$F(t)/G(t)\le F(x)/G(x)$ for $\xi x\le t\le x.$ Thus we can
construct $t_n\to\infty$ so that $F(t_n)/G(t_n)\to\infty$
and $F_{t_n}(x) \le G_{t_n}(x)\le x^{q_0}$ for $n^{-1}\le
x\le 1.$ Thus $\cal F_{\infty}$ contains a function $H$ with
$H(x)\le x^{q_0}$.  This contradicts the fact that
$q_0<p_F.$ Thus $F(t)/G(t)\to\infty$ and it follows easily
that since $\Cal F_{\infty}=\Cal F$ that $p_F=q_F$.  We can
invoke Theorem 5.3 to obtain that both $F$ and $G$ must be
elastic.

The case $p_1=q_F$ is similar.  In this case $G(t)/F(t)$ is
unbounded and we use the same argument as above to show that
$G(t)/F(t)\to \infty.$

We omit the case $p_F<p_G$; the reasoning is much the
same.\bull\enddemo

\demo{Example}It remains to construct an example of a
Calder\'on couple $(L_F[0,1],L_G[0,1])$ with $F$ and $G$
non-equivalent and $p_F=p_G<q_F=q_G.$ Such an example is a
counterexample to the previously mentioned conjecture of
Brudnyi [8].  Our construction depends on the following
lemma:

\proclaim{Lemma 6.10}Let $(Y_0,Y_1)$ be a Calder\'on couple
and let $X$ be a Banach space.  Then the pair $(X\oplus
Y_0,X\oplus Y_1)$ also forms a Calder\'on
couple.\endproclaim

\demo{Proof}We suppose the direct sums are $\ell_1-$sums.
Suppose $(x_0,y_0), (x_1,y_1)\in X\oplus(Y_0+Y_1)$ satisfies
$$K(t,(x_0,y_0),X\oplus Y_0,X\oplus Y_1)\le
K(t,(x_1,y_1),X\oplus Y_0,X\oplus Y_1).$$ Then we observe
that $$\|x_0\|_X \le \|x_1\|_X + K(1,y_1,Y_0,Y_1).$$ Thus
there is an operator $S:X\oplus (Y_0+Y_1)\to X$ with
$\|S_0\|\le 1$ and $S(x_1,y_1)=x_0.$ On the other hand:  $$
K(t,y_0,Y_0,Y_1) \le \min(1,t)\|x_1\| +K(t,y_1,Y_0,Y_1)$$
Now $(Y_0,Y_1)$ is Gagliardo complete ([13], Lemma 3) so by
K-divisibility ([4],[7], [15]) we can write $y_0=u+v$ where
$$ K(t,u,Y_0,Y_1)\le \gamma\min(1,t)\|x_1\| $$ and $$
K(t,v,Y_0,Y_1)\le \gamma K(t,y,Y_0,Y_1)$$ and $\gamma$ is an
absolute constant.  The former inequality implies that
$\max(\|u\|_{Y_0},\|u\|_{Y_1})\le \gamma\|x_1\|_X $ and
hence that there exists $S_1:X\to Y_0\cap Y_1$ with
$\|S_1\|\le \gamma $ and $S_1x_1=u.$ The latter inequality
yiels the existence of $S_2:Y_0+Y_1\to Y_0+Y_1$ with $S_2\in
\Cal A(Y_0,Y_1)$ and $S_2y_1=v.$ Let $S(x,y)= (S_0x,
S_1x+S_2y).$ Then $S$ is bounded on each $X\oplus Y_i$ and
maps $(x_1,y_1)$ to $(x_0,y_0).$ \enddemo

We now construct the example.  We suppose $q>p>1$; we set
$r=\frac12(p+q)$, $\alpha=p-1$ and $\beta=q-r.$ We next
define $a_1=1$ and then inductively $(b_n)_{n\ge
1},(c_n)_{n\ge 1},(d_n)_{n\ge 1}$ and $(a_n)_{n\ge 2}$ by
letting $b_n=2^na_n,$ $c_n=4b_n,$ $d_n=c_n+2n$ and
$a_{n+1}=4d_n.$

We then can construct an unbounded nonnegative Lipschitz
function $\phi:\bold R\to \bold R$ so that supp $\phi\subset
\cup_{n\ge 1}[a_n,b_n]$ and $|\phi'(x)| \le \alpha x^{-1}$
a.e.  (or equivalently $|\phi(x)-\phi(y)|\le |\log x-\log
y|$ for $x,y\ge 1.$) We then also define a nonnegative
Lipschitz function $\psi:\bold R\to\bold R$ with supp
$\psi\subset \cup_{n\ge 1}[c_n,d_n]$ by defining
$\psi(x)=\beta(x-c_n)$ for $c_n\le x\le c_n+n$ and
$\psi(x)=\beta(d_n-x)$ for $c_n+n\le x\le d_n.$ Finally we
put $F(x) = x^r\exp (\psi(\log x)$ and
$G(x)=x^r\exp(\psi(\log x)+\phi(\log x)).$

Now observe that $F$ and $G$ both satisfy the
$\Delta_2$-condition and both $F(x)/x$ and $G(x)/x$ are
increasing functions so that $F$ and $G$ are equivalent to
convex Orlicz functions.  We prefer to work directly with
$F$ and $G.$ We consider the pair $(E_F,E_G)$.  For $n<0$
let $\lambda_n$ be the unique solution of
$F(\lambda_n)=2^{-n}$ and let $\nu_n$ be the unique solution
of $G(\nu_n)=2^{-n}.$ We split $\bold Z_-$ into two disjoint
sets $J_0,J_1$ by setting $J_0=\{n:\log\lambda_n\in \cup_k
[c_k/2,2d_k]\}$ and $J_1=\bold Z_-\setminus J_0.$

We claim that on $\omega(J_0)$ the norms $\|\ \|_{L_F}$ and
$\|\ \|_{L_G}$ are equivalent.  In fact since $F\le G$ we
need only bound $\sum_{n\in J_0} 2^nG(\xi_n)$ subject to
$\sum_{n\in J_0}2^nF(\xi_n)=1.$ To do this observe that if
$n\in J_0$ and $0\le \xi_n\le \lambda_n$ then
$F(\xi_n)=G(\xi_n)$ unless $\log \xi_n< \log\lambda_n/2.$
Thus $$ \align \sum_{n\in J_0}G(\xi_n) &\le 1 + \sum_{n\in
J_0}2^nG(\sqrt\lambda_n) \\ &\le 1+ \sum_{n\in
J_0}\lambda_n^{-1/2} \endalign $$ and this establishes the
required estimate since $\lambda_n$ increases geometrically.

On $\omega(J_1)$ we claim both $\|\,\|_F$ and $\|\,\|_G$ are
equivalent to weighted $\ell_r-$norms and hence form a
Calder\'on couple by the result of Sparr [36].  Let us do
this for the case of $\|\,\|_G$ which we claim is equivalent
to $(\sum_{n\in J_1} |\xi_n/\nu_n|^r)^{1/r}.$ It suffices to
(a) bound $\sum_{n\in J_1}2^nG(\xi_n)$ subject to
$\sum_{n\in J_1}\xi_n^r\nu_n^{-r}=1$ and (b) conversely
bound $\sum_{n\in J_1}\xi_n^r \nu_n^{-r}$ subject to
$\sum_{n\in J_1}2^nG(\xi_n)=1.$

For (a) note that if $0\le \xi_n\le \nu_n$ then $$ |\log
G(\xi_n)-\log G(\nu_n) - r\log (\frac{\xi_n}{ \nu_n})|
\le\alpha \log(\frac{\log\nu_n}{\log\xi_n})$$ as long as
$\log\xi_n >\log \nu_n /2.$ Hence $$ G(\xi_n) \le
2^{\alpha}2^{-n}\nu_n^{-r}\xi_n^r+G(\sqrt\nu_n)$$ for $n\in
J_1.$ Thus $$ \sum_{n\in J_1}2^nG(\xi_n) \le 2^{\alpha}
+\sum_{n\in J_1}2^nG(\sqrt \nu_n)\le 2^{\alpha}+\sum_{n\in
J_1}\nu_n^{-1/2}$$ and this gives the required estimate.
(b) is similar.  The argument that $\|\ \|_F$ is equivalent
to $(\sum_{n\in J_1}|\xi_n|^r\lambda_n^{-r})^{1/r}$ is
slightly simpler and we omit it.

This completes the construction of the example.  It is clear
from Lemma 6.10 that $(E_F,E_G)$ and hence
$(L_F[0,1],L_G[0,1])$ is a Calder\'on couple with
$p_F=p_G=p$ and $q_F=q_G=q$ but that $F$ and $G$ are
non-equivalent.

We remark in closing that it is possible to find Orlicz
function spaces $L_F[0,1]$ so that if $(L_F,L_G)$ forms a
Calder\'on pair then $L_F=L_G.$ (We assume the
$\Delta_2-$condition for both $F$ and $G.$) We sketch the
details.  The argument of Theorem 6.9 can be used to
establish that if $F$ and $G$ are not equivalent at $\infty$
then there exists $p$ with $1\le p<\infty$ so that $x^p$ is
equivalent, for $0\le x\le 1,$ to a function of the form
$\lim F_{t_n}(x)$ where $t_n\to\infty$.  Now there are many
examples of functions $F$ which fail this property; for
example one can take the minimal Orlicz function:  $$ F(x) =
x^2\exp(\alpha\sum_{n=0}^{\infty}(1-\cos(2\pi(\log
t)/2^n)).$$ See [20].

\subheading{Acknowledgements} The author wishes to thank
Michael Cwikel for introducing him to the problems studied
in this paper and for many stimulating discussions and
helpful comments.  He would also like to thank E. Pustylnik
for some valuable comments on an earlier draft of the paper.

\subheading{References}

\item{1.}J.  Arazy and M. Cwikel, A new characterization of
the interpolation spaces betwen $L_p$ and $L_q$, Math.
Scand. 55 (1984) 253-270.

\item{2.}S.F.  Bellenot, The Banach spaces of Maurey and
Rosenthal and totally incomparable bases, J. Functional
Analysis 95 (1991) 96-105.

\item{3.}S.F.  Bellenot, R. Haydon and E. Odell,
Quasi-reflexive and tree spaces constructed in the spirit of
R.C.  James, Contemporary Math. 85 (1987) 19-43.

\item{4.}C.  Bennett and R. Sharpley, {\it Interpolation of
Operators,} Academic Press, New York 1988.

\item{5.}J.  Bergh and J. L\"ofstr\"om, {\it Interpolation
spaces.  An Introduction,} Springer, Berlin 1976.

\item{6.}N.H.  Bingham, C.M.  Goldie and J.L.  Teugels, {\it
Regular Variation,} Cambridge University Press, Cambridge
1987.

\item{7.}Y.  Brudnyi and N. Krugljak, Real interpolation
functors, Soviet Math.  Doklady 23 (1981) 5-8.

\item{8.}Y.  Brudnyi and N. Krugljak, {\it Interpolation
functors and interpolation spaces} North Holland 1991.

\item{9.}A.P.  Calder\'on, Spaces between $L_1$ and
$L_{\infty}$ and the theorems of Marcinkiewicz, Studia Math.
26 (1966) 273-299.

\item{10.}P.G.  Casazza, W.B.  Johnson and L. Tzafriri, On
Tsirelson's space, Israel J. Math. 47 (1984) 81-98.

\item{11.}P.G.  Casazza and B.L.  Lin, On symmetric basic
sequences in Lorentz sequence spaces, II, Israel J. Math. 17
(1974) 191-218.

\item{12.}P.G.  Casazza and T.J.  Shura, {\it Tsirelson's
space,} Springer Lecture Notes 1363, Berlin 1989

\item{13.}M.  Cwikel, Monotonicity properties of
interpolation spaces, Ark.  Mat. 14 (1976) 213-236.

\item{14.}M.  Cwikel, Monotonicity properties of
interpolation spaces II, Ark.  Mat. 19 (1981) 123-136.

\item{15.}M.  Cwikel, K-divisibility of the K-functional and
Calder\'on couples, Ark.  Mat. 22 (1984) 39-62.

\item{16.}M.  Cwikel and P. Nilsson, On Calder\'on-Mityagin
couples of Banach lattices, Proc.  Conf.  Constructive
Theory of Functions, Varna, 1984, Bulgarian Acad.  Sciences
1984, 232-236.

\item{17.}M.  Cwikel and P. Nilsson, Interpolation of
Marcinkiewicz spaces, Math.  Scand. 56 (1985) 29-42.

\item{18.}M.  Cwikel and P. Nilsson, Interpolation of
weighted Banach lattices, Memoirs Amer.  Math.  Soc. to
appear.

\item{19.}G.H.  Hardy, J.E.  Littlewood and G. Polya, {\it
Inequalities,} Cambridge University Press, 1934.

\item{20.} F.L.  Hernandez and B. Rodriguez-Salinas, On
$\ell^p$-complemented copies in Orlicz spaces II, Israel J.
Math. 68 (1989) 27-55.

\item{21.}P.W.  Jones, On interpolation between $H_1$ and
$H_{\infty}$, 143-151 in {\it Interpolation spaces and
allied topics in analysis, Proceedings, Lund Conference
1983,} M. Cwikel and J. Peetre, editors, Springer Lecture
Notes 1070, 1984.

\item{22.}J.L.  Krivine, Sous espaces de dimension finite
des espaces de Banach reticul\'es, Ann.  Math. 104 (1976)
1-29.

\item{23.}J.  Lindenstrauss and L. Tzafriri, On the
complemented subspaces problem, Israel J. Math. 9 (1971)
263-269.

\item{24.}J.  Lindenstrauss and L. Tzafriri, {\it Classical
Banach Spaces I} Springer, Berlin 1977.

\item{25.}J.  Lindenstrauss and L. Tzafriri, {\it Classical
Banach Spaces II} Springer, Berlin 1979.

\item{26.}G.G.  Lorentz, Relations between function spaces,
Proc.  Amer.  Math.  Soc. 12 (1961) 127-132.

\item{27.}G.G.  Lorentz and T. Shimogaki, Interpolation
theorems for the pairs of spaces $(L_1,L_p)$ and
$(L_p,L_{\infty})$, Trans.  Amer.  Math.  Soc. 159 (1971)
207-222.

\item{28.}L.  Maligranda, On Orlicz results in interpolation
theory, Proc.  Orlicz Memorial Conference, Univ. of
Mississippi, 1990.

\item{29.}L.  Maligranda and V.I.  Ovchinnikov, On
interpolation between $L_1+L_{\infty}$ and $L_1\cap
L_{\infty}$, J. Functional Analysis 107 (1992) 342-351.

\item{30.}C.  Merucci, Interpolation r\'eelle avec fonction
param\`etre:  r\'eit\'eration et applications aux espaces
$\Lambda^p(\phi),\ (0<p\le +\infty)$, C.R.  Acad.  Sci.
(Paris) I, 295 (1982) 427-430.

\item{31.}C.  Merucci, Applications of interpolation with a
function parameter to Lorentz, Sob\-ol\-ev and Besov spaces,
183-201 in {\it Interpolation spaces and allied topics in
analysis, Proceedings, Lund Conference 1983,} M. Cwikel and
J. Peetre, editors, Springer Lecture Notes 1070, 1984.

\item{32.}B.S.  Mityagin, An interpolation theorem for
modular spaces, Mat. Sbornik 66 (1965) 473-482.

\item{33.}S.J.  Montgomery-Smith, Comparison of
Orlicz-Lorentz spaces, Studia Math. to appear.

\item{34.}V.I.  Ovchinnikov, On the estimates of
interpolation orbits, Mat. Sb. 115 (1981) 642-652 (= Math.
USSR Sbornik 43 (1982) 573-583.)

\item{35.}A.A.  Sedaev and E.M.  Semenov, On the possibility
of describing interpolation spaces in terms of Peetre's
K-method, Optimizaciya 4 (1971) 98-114.

\item{36.}G.  Sparr, Interpolation of weighted $L_p-$spaces,
Studia Math. 62 (1978) 229-271.

\item{37.}B.S.  Tsirelson, Not every Banach space contains
an embedding of $\ell_p$ or $c_0$, Functional Anal.  Appl. 8
(1974) 138-141.

\item{38.}Q.  Xu, Notes on interpolation of Hardy spaces,
preprint, 1991.

\item{39.}M.  Zippin, On perfectly homogeneous bases in
Banach spaces, Israel J. Math. 4 (1966) 265-272.

\bye